\documentclass[11pt,leqno]{article}

\usepackage{amssymb,amsmath,amsthm,mathrsfs}
\usepackage{graphicx}
\usepackage{times}

\leftmargini=\baselineskip

\newtheoremstyle{localthm}
	{5pt} % space above
	{5pt} % space below
	{\sl} % Body font
	{} % Indent amount
	{\bf} % Theorem head font
	{{\rm.}} % Punctuation after theorem head
	{.7em} % Space after theorem head
	{} % Theorem head spec ?

\theoremstyle{localthm}
\newtheorem{Theorem}{Theorem}[section]
\newtheorem{Corollary}[Theorem]{Corollary}

\newtheorem{Lemma}[Theorem]{Lemma}

\def\GGiso{\GG^{}_\uparrow}
\def\GGconv{\GG^{}_{\rm conv}}
\def\Dk{\nabla^k}

\def\nin{\noindent}

\def\be{\begin{equation}}
\def\ee{\end{equation}}
\def\bea{\begin{eqnarray*}}
\def\eea{\end{eqnarray*}}
\def\bean{\begin{eqnarray}}
\def\eean{\end{eqnarray}}

\def\barr{\begin{array}}
\def\earr{\end{array}}

\def\CC{{\mathcal C}}
\def\FF{{\mathcal F}}
\def\GG{{\mathcal G}}
\def\HH{{\mathcal H}}
\def\LL{{\mathcal L}}
\def\NN{{\mathcal N}}
\def\TT{{\mathcal T}}

\def\Bl{\Bigl}
\def\Br{\Bigr}

\def\R{\mathbb{R}}

\def\Ex{\mathrm{I\!\!E}}
\def\Pr{\mathrm{I\!\!P}}
\def\Var{\mathrm{Var}}

\begin{document}
%===============

\addtolength{\baselineskip}{+.4\baselineskip}
\addtolength{\parskip}{+.2\baselineskip}

\title{\bf Optimal Confidence Bands \\ for Shape-Restricted Curves}
\author{Lutz D\"umbgen\\\\
	Department of Mathematical Statistics and Actuarial Science\\
	University of Bern\\
	Sidlerstrasse 5, CH-3012 Bern, Switzerland\\\\
	E-mail: duembgen@stat.unibe.ch, URL: www.imsv.unibe.ch}
\date{August 2002, updated December 2013}
\maketitle

\begin{center}
This paper has been published in \textsl{Bernoulli \textbf{9}}, 423--449 (2003).\\
The present version corrects two typos in the published version\\
and contains updated references.
\end{center}

\vfill

\centerline{\bf Abstract}
Let $Y$ be a stochastic process on $[0,1]$ satisfying $dY(t) = n^{1/2} f(t) dt + dW(t)$, where $n \ge 1$ is a given scale parameter (``sample size''), $W$ is standard Brownian motion and $f$ is an unknown function. Utilizing suitable multiscale tests we construct confidence bands for $f$ with guaranteed given coverage probability, assuming that $f$ is isotonic or convex. These confidence bands are computationally feasible and shown to be asymptotically sharp optimal in an appropriate sense.

\vfill

\nin{\bf Running title.} Confidence Bands for Shape-Restricted Curves

\nin{\bf Keywords and phrases.} adaptivity, concave, convex, isotonic, kernel estimator, local smoothness, minimax bounds, multiscale testing

\newpage

\section{Introduction}
\label{Introduction} 

Nonparametric statistical models often involve some unknown function $f$ defined on a real interval $J$. For instance $f$ might be the probability density of some distribution or a regression function. Nonparametric point estimators for such a curve $f$ are abundant. The available methods are based on kernels, splines, local polynomials, or orthogonal series, including wavelets; see Hart~(1997) and references cited therein. In order to quantify the precision of estimation, one often wants to replace a point estimator with a confidence band $(\hat\ell, \hat u)$ for $f$. The latter consists of two functions $\hat\ell = \hat\ell(\cdot,\mbox{data})$ and $\hat u = \hat u(\cdot,\mbox{data})$ on $J$ with values in $[-\infty,\infty]$ such that, hopefully, $\hat\ell \le f \le \hat u$ pointwise. More precisely, one is aiming at a confidence band such that
\be
	\Pr\{\hat\ell \le f \le \hat u\} \ \ge \ 1 - \alpha	
	\label{Confidence0}
\ee
for a given level $\alpha \in \left]0,1\right[$, while $\hat \ell$ and $\hat u$ should be as close to each other as possible.

Unfortunately, curve estimation is an ill-posed problem, and usually there are no nontrivial bands $(\hat\ell,\hat u)$ satisfying (\ref{Confidence0}) for arbitrary $f$; see Donoho~(1988). Therefore one has to impose some additional restrictions on $f$. One possibility are smoothness constraints on $f$, for instance an upper bound on a certain derivative of $f$. Under such restrictions, (\ref{Confidence0}) can be achieved approximately for large sample sizes; see for example Bickel and Rosenblatt~(1973), Knafl et al.~(1985), Hall and Titterington~(1988), H\"ardle and Marron~(1991), Eubank and Speckman~(1993), Fan and Zhang~(2000), and the references cited therein.

A problem with the aforementioned methods is that smoothness constraints are hard to justify in practical situations. More precisely, even if the underlying curve $f$ is infinitely often differentiable, the actual coverage probabilities of the confidence bands mentioned above depend on quantitative properties of certain derivatives of $f$ which are difficult to obtain from the data.

In many applications qualitative assumptions about $f$ such as monotonicity, unimodality or concavity/convexity are plausible. One example are growth curves in medicine, e.g.~where $f(x)$ is the mean body height of newborns at age $x$. Here isotonicity of $f$ is a plausible assumption. Another example are so-called Engel curves in econometrics, where $f(x)$ is the mean expenditure for certain consumer goods of households with annual income $x$. Here one expects $f$ to be isotonic and sometimes concave as well. Under such qualitative assumptions it is possible to construct $(1 - \alpha)$--confidence sets for $f$ based on certain goodness-of-fit tests without relying on asymptotic arguments. Examples for such procedures can be found in Davies~(1995), Hengartner and Stark~(1995) and D\"umbgen~(1998). In particular, these papers present confidence bands $(\hat\ell,\hat u)$ for $f$ such that
\be
	\Pr\{\hat\ell \le f \le \hat u\} \ \ge \ 1 - \alpha	
	\quad\mbox{whenever } f \in \FF .
	\label{Confidence}
\ee
Here $\FF$ denotes the specified class of functions. Given a suitable distance measure $D(\cdot,\cdot)$ for functions, the goal is to find a band $(\hat\ell,\hat u)$ satisfying (\ref{Confidence}) such that either $D(\hat u, \hat \ell)$ or $D(\hat\ell, f)$ and $D(\hat u,f)$ are as small as possible. The phrase ``as small as possible'' can be interpreted in the sense of optimal rates of convergence to zero as the sample size $n$ tends to infinity. The papers of Hengartner and Stark~(1995) and D\"umbgen~(1998) contain such optimality results. 

In the present paper we investigate optimality of confidence bands in more detail. In addition to optimal rates of convergence we obtain optimal constants and discuss the impact of local smoothness properties of $f$. Compared to the general confidence sets of D\"umbgen~(1998), the methods developed here are more stringent and computationally simpler. They are based on multiscale tests as developed by D\"umbgen and Spokoiny~(2001), who considered tests of qualitative assumptions rather than confidence bands. For further results on testing in nonparametric curve estimation see Hart~(1997), Fan et al.~(2001), and the references cited there.

%===================================
\section{Basic setting and overview}
%===================================

For mathematical convenience we focus on a continuous white noise model: Suppose that one observes a stochastic process $Y$ on the unit interval $[0,1]$, where
$$
    Y(t) \ = \ n^{1/2} \int_0^t f(x) \, dx + W(t) .
$$
Here $f$ is an unknown function in $L^2[0,1]$, $n \ge 1$ is a given scale parameter (``sample size''), and $W$ is standard Brownian motion. In this context the bounding functions $\hat\ell, \hat u$ are defined on $[0,1]$, but for notational convenience the function $f$ is tacitly assumed to be defined on the whole real line with values in $[-\infty,\infty]$. From now on we assume that
$$
	f \ \in \ \GG \cap L^2[0,1] ,
$$
where $\GG$ denotes one of the following two function classes:
\bea
	\GGiso & := & \Bl\{ \mbox{non-decreasing functions } g : \R \to [-\infty,\infty] \Br\} , \\
	\GGconv & := & \Bl\{ \mbox{convex functions } g : \R \to \left]-\infty,\infty\right] \Br\} .
\eea

The paper is organized as follows. In Section~\ref{GGiso Levy} we treat the case $\GG = \GGiso$ and measure the quality of a confidence band $(\hat\ell, \hat u)$ by quantities related to the Levy distance $d_{\rm L}(\hat\ell,\hat u)$. Generally,
$$
	d_{\rm L}(g,h) \ := \ \inf \Bl\{ \epsilon > 0 : g \le h(\cdot + \epsilon) + \epsilon \mbox{ and } 
		h \le g(\cdot + \epsilon) + \epsilon \mbox{ on } [0,1-\epsilon] \Br\}
$$
for isotonic functions $g,h : [0,1] \to [-\infty,\infty]$. It turns out that a confidence band which is based on a suitable multiscale test as introduced by D\"umbgen and Spokoiny~(2001) is asymptotically optimal in a strong sense. Throughout this paper asymptotic statements refer to $n\to\infty$, unless stated otherwise. 

In Section~\ref{GG smooth} we treat both classes $\GGiso$ and $\GGconv$ simultaneously. We discuss the construction of confidence bands $(\hat\ell,\hat u)$ satisfying (\ref{Confidence}) such that $D(\hat\ell,f)$ and $D(f,\hat u)$ are as small as possible whenever $f$ satisfies some additional smoothness constraints. Here $D(g,h)$ is a distance measure of the form
$$
	D(g,h) \ := \ \sup_{x \in [0,1]} \, w(x,f) (h(x) - g(x))
$$
for some weight function $w(\cdot,f) \ge 0$ reflecting local smoothness properties of $f$. Again it turns out that suitable multiscale procedures yield nearly optimal procedures without additional prior information on $f$.

In Section~\ref{Examples} we present some numerical examples for the procedures of Section~\ref{GG smooth}. The proofs are deferred to Sections~\ref{Proofs}, \ref{Decision Theory} and \ref{Optimization}. In particular, Section~\ref{Decision Theory} contains a new minimax bound for confidence rectangles in a gaussian shift model, which may be of independent interest.

As for the white noise model, the results of Brown and Low~(1996), Nussbaum~(1996) and Grama
and Nussbaum~(1998) on asymptotic equivalence can be used to transfer the lower bounds of the present paper to other models. Moreover, one can mimick the confidence bands developed here in traditional regression models under minimal assumptions; see D\"umbgen and Johns~(2004) and D\"umbgen~(2007).

%===========================================================
\section{Optimality for isotonic functions in terms of 
	L\'{e}vy type distances}
\label{GGiso Levy}
%===========================================================

\begin{figure}[h]
\centering
\includegraphics[width=7cm,height=9cm]{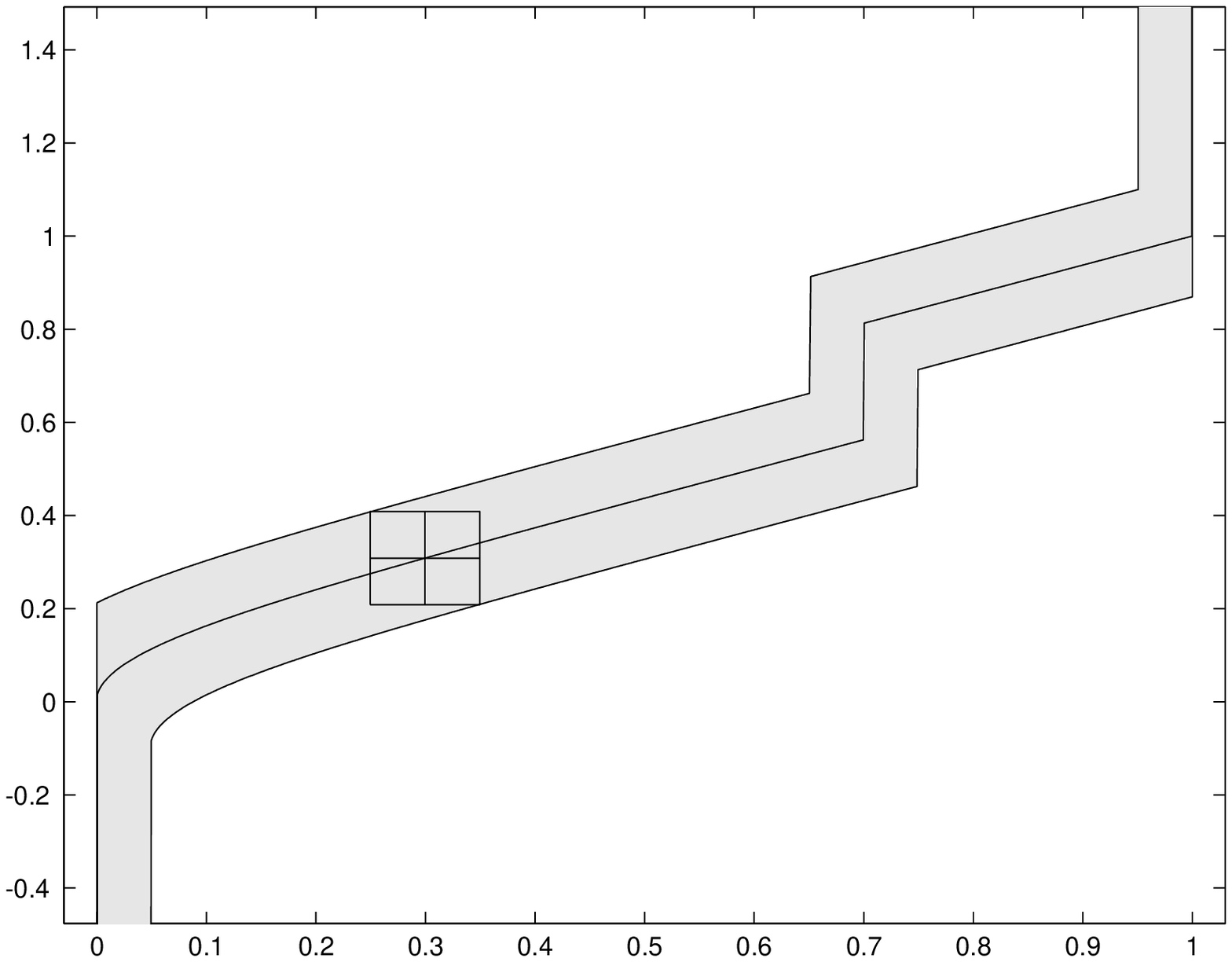}
\hfill 
\includegraphics[width=7cm,height=9cm]{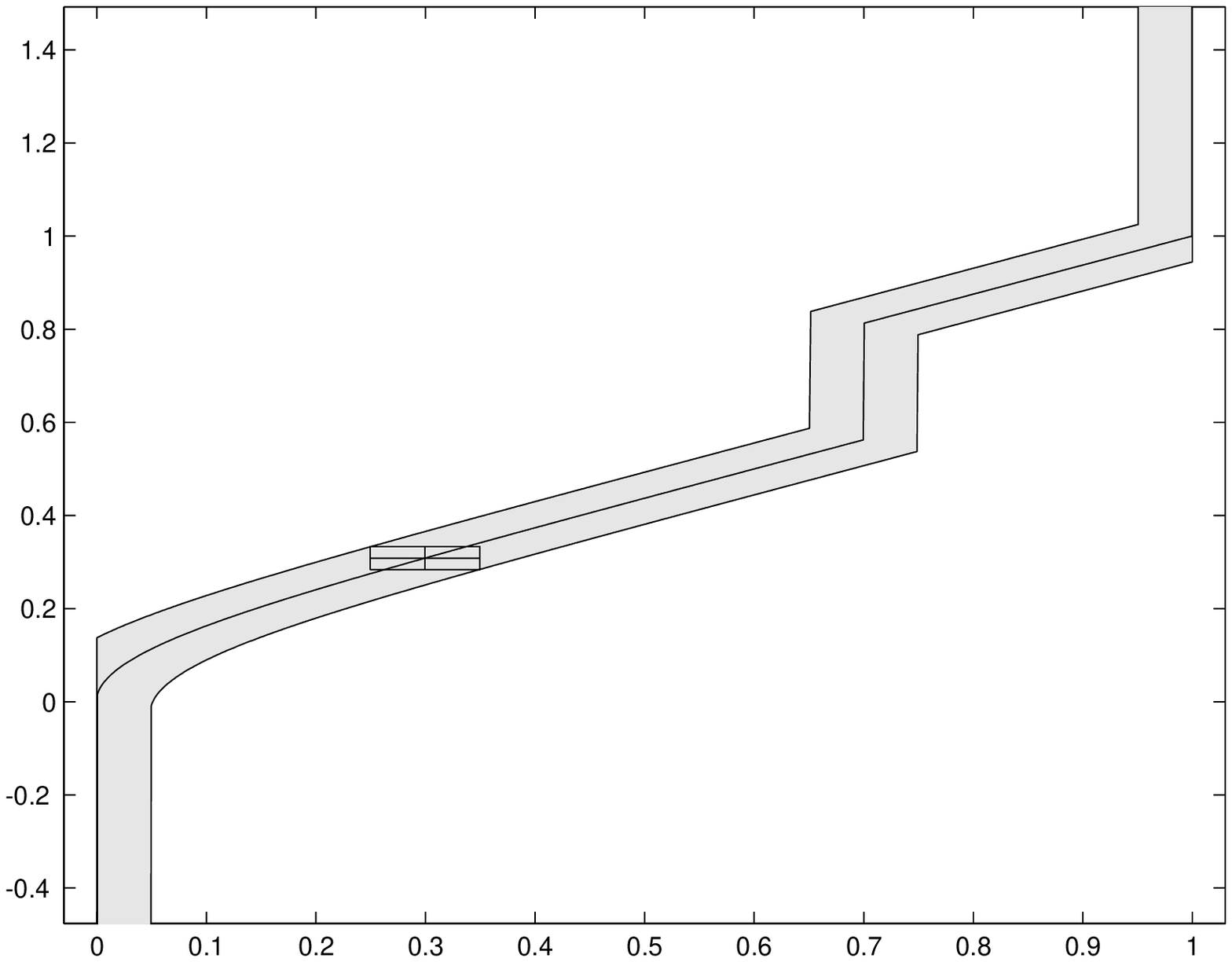}
\caption{Two $D_{0.05}(\cdot,\cdot)$--neighborhoods of some function $g$.}
\label{Levy-Band}
\end{figure}

In this section we consider the class $\GGiso$. For isotonic functions $g,h : [0,1] \to [-\infty,\infty]$ and $\epsilon > 0$ let
$$
	D_\epsilon(g,h) 
	\ := \ \inf \Bl\{ \lambda \ge 0 : 
		g \le h(\cdot + \epsilon) + \lambda \mbox{ and } h \le g(\cdot + \epsilon) + \lambda 
		\mbox{ on } [0,1-\epsilon] \Br\} .
$$
Then the L\'{e}vy distance $d_{\rm L}(g,h)$ is the infimum of all $\epsilon > 0$ such that $D_\epsilon(g,h) \le \epsilon$. We use these functionals $D_\epsilon(\cdot,\cdot)$ in order to quantify differences between isotonic functions. Figure~\ref{Levy-Band} depicts one such function $g$, and the shaded areas represent the set of all functions $h$ with $D_{0.05}(g,h) \le 0.1$ and $D_{0.05}(g,h) \le 0.025$, respectively.

The next theorem provides lower bounds for $D_\epsilon(\hat\ell, \hat u)$, $0 < \epsilon \le 1$. Here and throughout the sequel the dependence of probabilities, expectations and distributions on the functional parameter $f$ is sometimes indicated by a subscript $f$.

\begin{Theorem}	\label{Lower Bounds GGiso}
There exists a universal function $b$ on $\left]0,1\right]$ with $\lim_{\epsilon \downarrow 0} b(\epsilon) = 0$ such that
$$
	\inf_{f \in \GGiso \cap L^2[0,1]} \, 
	\Pr_f \left\{ \hat\ell \le f \le \hat u \mbox{ and } 
		D_\epsilon(\hat\ell, \hat u) 
			< { (8 \log(e/\epsilon))^{1/2} - b(\epsilon) \over (n\epsilon)^{1/2}} \right\} 
	\ \le \ b(\epsilon)
$$
for any confidence band $(\hat\ell, \hat u)$ and arbitrary $\epsilon \in \left]0,1\right]$.
\end{Theorem}

Theorem~\ref{Lower Bounds GGiso} entails a lower bound for $d_{\rm L}(\hat\ell,\hat u)$. For let $\epsilon = \epsilon_n := c \, (\log(n)/n)^{1/3} - \delta n^{-1/3}$ with any fixed $c, \delta > 0$. Then one can show that for sufficiently large $n$,
$$
	{(8 \log(e/\epsilon))^{1/2} - b(\epsilon) \over (n\epsilon)^{1/2}} 
	\ = \ \Bl( {8 \over 3c} \Br)^{1/2} \Bl( {\log n\over n} \Br)^{1/3} + o(n^{-1/3})
	\ \ge \ \epsilon ,
$$
provided that $c$ equals $(8/3)^{1/3} \approx 1.387$.

\begin{Corollary}
For each $n \ge 1$ there exists a universal constant $\beta_n$ such that $\beta_n\to 0$ and
$$
	\inf_{f \in \GGiso \cap L^2[0,1]} \, 
	\Pr_f \left\{ \hat\ell \le f \le \hat u \mbox{ and } 
		d_{\rm L}(\hat\ell, \hat u) 
			< \Bl( {8\over 3} \Br)^{1/3} \Bl( {\log n\over n} \Br)^{1/3} - \beta_n n^{-1/3} \right\} 
	\ \le \ \beta_n
$$
for any confidence band $(\hat\ell, \hat u)$.	\hfill	$\Box$
\end{Corollary}

It is possible to get close to these lower bounds for $D_\epsilon(\hat\ell, \hat u)$ {\sl simultaneously} for all $\epsilon \in \left]0,1\right]$ while (\ref{Confidence}) is satisfied. For let $\kappa_\alpha$ be a real number such that
$$
	\Pr \left\{ {|W(t) - W(s)| \over (t - s)^{1/2}} \le \Gamma(t - s) + \kappa_\alpha 
		\mbox{ for } 0 \le s < t \le 1 \right\} 
	\ \le \ \alpha ,
$$
where
$$
	\Gamma(u) \ := \ (2 \log(e/u))^{1/2}	\quad\mbox{for } 0 < u \le 1 .
$$
The existence of such a critical value $\kappa_\alpha$ follows from D\"umbgen and Spokoiny~(2001, Theorem~2.1). With the local averages
$$
	F_f(s,t) \ := \ {1 \over t - s} \int_s^t f(x) \, dx
$$
of $f$ and their natural estimators
$$
	\hat F(s,t)	\ := \ {Y(t) - Y(s) \over n^{1/2}(t - s)}
$$
it follows that
$$
	\Pr_f \left\{ \Bl| \hat F(s,t) - F_f(s,t) \Br| 
		\le {\Gamma(t - s) + \kappa_\alpha \over (n(t - s))^{1/2}} 
		\mbox{ \ for } 0 \le s < t \le 1 \right\} 
	\ \ge \ 1 - \alpha .
$$
But for $0 \le s < t \le 1$,
$$
	f(s) \ \le \ F_f(s,t) \ \le \ f(t)	
	\quad\mbox{whenever $f \in \GGiso$} .
$$
This implies the first assertion of the following theorem.

\begin{Theorem}	\label{Upper Bounds GGiso}
With the critical value $\kappa_\alpha$ above let
\bea
	\hat\ell(x) 
	& := & \sup_{0 \le s < t \le x} \, 
		\Bl( \hat F(s,t) - {\Gamma(t - s) + \kappa_\alpha \over \sqrt{n(t - s)}} \Br) , \\
	\hat u(x) 
	& := & \inf_{x \le s < t \le 1} \, 
		\Bl( \hat F(s,t) + {\Gamma(t - s) + \kappa_\alpha \over \sqrt{n(t - s)}} \Br) .
\eea
This defines a confidence band $(\hat\ell, \hat u)$ for $f$ satisfying (\ref{Confidence}) with $\FF = \GGiso \cap L^2[0,1]$. Moreover, in case of $\hat\ell \le \hat u$,
\bea
	D_\epsilon(\hat\ell, \hat u) 
	& \le & {(8 \log(e/\epsilon))^{1/2} + 2 \kappa_\alpha \over (n\epsilon)^{1/2}}	
		\quad\mbox{for } 0 < \epsilon \le 1 , \\
	d_{\rm L}(\hat\ell, \hat u) 
	& \le & \Bl( {8\over 3} \Br)^{1/3} \Bl( {\log n\over n} \Br)^{1/3} + o(n^{-1/3}) .
\eea
\end{Theorem}

\noindent
{\bf Proof.} The preceding upper bound for $D_\epsilon(\hat\ell, \hat u)$ follows from the fact that for any $x \in [0,1-\epsilon]$,
\bea
	\hat u(x) - \hat\ell(x+\epsilon) 
	& \le & \Bl( \hat F(x,x+\epsilon)
			+ {\Gamma(\epsilon) + \kappa_\alpha \over (n\epsilon)^{1/2}} \Br) 
		- \Bl( \hat F(x,x+\epsilon)
			- {\Gamma(\epsilon) + \kappa_\alpha \over (n\epsilon)^{1/2}} \Br) \\
	& = & {2 \Gamma(\epsilon) + 2 \kappa_\alpha \over (n\epsilon)^{1/2}} \\
	& = & {(8 \log(e/\epsilon))^{1/2} + 2 \kappa_\alpha \over (n\epsilon)^{1/2}} .
\eea
Letting $\epsilon = \epsilon_n = (8/3)^{1/3} (\log(n)/n)^{1/3}$ yields the upper bound for $d_{\rm L}(\hat\ell, \hat u)$.	\hfill	$\Box$

%===============================================
\section{Bands for potentially smooth functions}
\label{GG smooth}
%===============================================

A possible criticism of the preceding results is the fact that the minimax bounds are attained at special step functions. On the other hand one often expects the underlying curve $f$ to be smooth in some vague sense. Therefore we aim now at confidence bands satisfying (\ref{Confidence}) with $\FF = \GG \cap L^2[0,1]$, which are as small as possible whenever $f$ satisfies some additional smoothness conditions. Throughout $\GG$ stands for $\GGiso$ or $\GGconv$.

In the sequel let $\langle g,h \rangle := \int_{-\infty}^\infty g(x)h(x) \, dx$ and $\|g\| := \langle g,g\rangle^{1/2}$ for measurable functions $g,h$ on the real line such that these integrals are defined. The confidence bands to be presented here can be described either in terms of kernel estimators for $f$ or in terms of tests. Both viewpoints have their own merits.

%=====================================
\subsection{Kernel estimators for $f$}
%=====================================

Let $\psi$ be some kernel function in $L^2(\R)$. For technical reasons we assume that $\psi$ satisfies the following three regularity conditions:
\be
	\left\{\barr{l} 
	\mbox{$\psi$ has bounded total variation}; \\
	\mbox{$\psi$ is supported by $[-a,b]$, where $a,b \ge 0$}; \\
	\langle 1,\psi \rangle \ > \ 0 .
	\earr\right.
	\label{Regularity}
\ee
For any bandwidth $h > 0$ and location parameter $t \in \R$ let
$$
	\psi_{h,t}(x) \ := \ \psi \Bl( {x - t \over h} \Br) .
$$
Then $\langle g, \psi_{h,t}\rangle = h \, \langle g(t + h \, \cdot), \psi\rangle$ and $\|\psi_{h,t}\| = h^{1/2} \|\psi\|$. A kernel estimator for $f(t)$ with kernel function $\psi$ and bandwidth $h$ is given by
$$
	\hat f_h(t) \ := \ {\psi Y(h,t) \over n^{1/2} h \, \langle 1,\psi\rangle} ,
$$
where
$$
	\psi Y(h,t) \ := \ \int_0^1 \psi_{h,t}(x) \, dY(x) .
$$
From now on suppose that $ah \le t \le 1 - bh$. Then $\psi_{h,t}$ is supported by $[0,1]$ and one may write
\bea
	\Ex \hat f_h(t) 
	& = & {\langle f, \psi_{t,h}\rangle \over h \, \langle 1,\psi\rangle} 
		\ = \ {\langle f(t + h\,\cdot), \psi\rangle \over \langle 1,\psi\rangle} , \\
	\Var(\hat f_h(t))
	& = & {\|\psi_{t,h}\|^2 \rangle \over n h^2 \langle 1,\psi\rangle^2}
		\ = \ {\|\psi\|^2 \over nh \, \langle 1,\psi\rangle^2} .
\eea

The random fluctuations of these kernel estimators can be bounded uniformly in $h > 0$. For that purpose we define the multiscale statistic
\bea
	T(\pm\psi)
	& := & \sup_{h > 0} \, \sup_{t \in [ah,1-bh]} 
		\Bl( {\pm \psi W(h,t) \over h^{1/2} \|\psi\|} - \Gamma((a+b)h) \Br) \\
	& = & \sup_{h > 0} \, \sup_{t \in [ah,1-bh]} 
		\Bl( \pm \, {\hat f_h(t) - \Ex \hat f_h(t) \over \Var(\hat f_h(t))^{1/2}}
			- \Gamma((a+b)h) \Br) ,
\eea
similarly as in D\"umbgen and Spokoiny~(2001). It follows from Theorem~2.1 in the latter paper, that $0 \le T(\pm\psi) < \infty$ almost surely. In particular, $|\hat f_h(t) - \Ex \hat f_h(t)| \le (nh)^{-1/2} \log(e/h)^{1/2} O_p(1)$, uniformly in $h > 0$ and $ah \le t \le 1 - bh$.

It is well-known that kernel estimators are biased in general. But our shape restrictions may be used to construct two kernel estimators whose bias is always non-positive or non-negative, respectively. Precisely, let $\psi^{(\ell)}$ and $\psi^{(u)}$ be two kernel functions satisfying (\ref{Regularity}) with respective supports $[-a^{(\ell)}, b^{(\ell)}]$ and $[-a^{(u)}, b^{(u)}]$. In addition suppose that
\bean
	\langle g, \psi^{(\ell)}\rangle
	& \le & g(0) \langle 1,\psi^{(\ell)}\rangle
		\quad\mbox{for all } g \in \GG \cap L^2[-a^{(\ell)},b^{(\ell)}] ,
	\label{OrthogonalityL} \\
	\langle g, \psi^{(u)}\rangle
	& \ge & g(0) \langle 1,\psi^{(u)}\rangle
		\quad\mbox{for all } g \in \GG \cap L^2[-a^{(u)},b^{(u)}] .
	\label{OrthogonalityU}
\eean
These inequalities imply that the corresponding kernel estimators satisfy the inequalities $\Ex \hat f_h^{(\ell)}(t) \le f(t) \le \Ex \hat f_h^{(u)}(t)$, and the definition of $T(\pm\psi)$ yields that
\bean
	f(t)
	& \ge & \hat f^{(\ell)}_h(t) 
		- {\|\psi^{(\ell)}\| \Bl( \Gamma(d^{(\ell)} h) + T(\psi^{(\ell)}) \Br)
			\over \langle 1,\psi^{(\ell)}\rangle (nh)^{1/2}} ,
	\label{Lower Bound} \\
	f(t)
	& \le & \hat f^{(u)}_h(t) 
		+ {\|\psi^{(u)}\| \Bl( \Gamma(d^{(u)} h) + T(- \psi^{(u)}) \Br)
			\over \langle 1,\psi^{(u)}\rangle (nh)^{1/2}} .
	\label{Upper Bound}
\eean
Here $d^{(z)} := a^{(z)} + b^{(z)}$. Now let $\kappa_\alpha$ be the $(1-\alpha)$--quantile of the combined statistic $T^* := \max \Bl( T(\psi^{(\ell)}), T(-\psi^{(u)}) \Br)$, i.e. the smallest real number such that $\Pr\{T^* \le \kappa_\alpha\} \ge 1 - \alpha$. Then
\bea
	\hat\ell(t) 
	& := & \sup_{h > 0 \ : \ t \in \left[a^{(\ell)}h,1-b^{(\ell)}h\right]} 
		\left( \hat f_h^{(\ell)}(t) 
		- {\|\psi^{(\ell)}\| (\Gamma(d^{(\ell)} h) + \kappa_\alpha)
			\over \langle 1,\psi^{(\ell)}\rangle (nh)^{1/2}} \right) , \\
	\hat u(t) 
	& := & \inf_{h > 0 \ : \ t \in \left[a^{(u)}h,1-b^{(u)}h\right]} 
		\left( \hat f_h ^{(u)}(t) 
		+ {\|\psi^{(u)}\| (\Gamma(d^{(u)} h) + \kappa_\alpha)
			\over \langle 1,\psi^{(u)}\rangle (nh)^{1/2}} \right)
\eea
defines a confidence band $(\hat\ell,\hat u)$ for $f$ satisfying (\ref{Confidence}).

Equality holds in (\ref{Confidence}) if $\GG = \GGiso$ and $f$ is constant, or if $\GG = \GGconv$ and $f$ is linear, provided that $\kappa_\alpha > 0$. For then it follows from (\ref{OrthogonalityL}) and (\ref{OrthogonalityU}) with $g(x) = \pm 1$ or $g(x) = \pm x$ that the kernel estimators are unbiased. Thus $\hat\ell \le f \le \hat u$ is equivalent to $T^*  > \kappa_\alpha$. Moreover, using general theory for gaussian measures on Banach spaces one can show that the distribution of $T^*$ is continuous on $\left]0,\infty\right[$.

Sufficient conditions for requirements~(\ref{OrthogonalityL}) and (\ref{OrthogonalityU}) in general are provided by Lemma~\ref{Ortho} in Section~\ref{Optimization}. The confidence band presented in Section~\ref{GGiso Levy} is a special case of the one derived here, if we define $\psi^{(\ell)}(x) := 1\{x \in [-1,0]\}$ and $\psi^{(u)}(x) := 1\{x \in [0,1]\}$ and apply postprocessing as described below.

%===================================================
\subsection{Postprocessing of confidence bands}
%===================================================

Any confidence band $(\hat\ell,\hat u)$ for $f$ can be enhanced, if we replace $\hat\ell(x)$ and $\hat u(x)$ with
$$
	\hat{\hat\ell}(x) \ := \ \inf \Bl\{ g(x) : g \in \GG, \hat\ell \le g \le \hat u \Br\}	
 \quad\mbox{and}\quad	
	\hat{\hat u}(x)   \ := \ \sup \Bl\{ g(x) : g \in \GG, \hat\ell \le g \le \hat u \Br\} ,
$$
respectively. Here we assume tacitly that the set $\{g \in \GG : \hat\ell \le g \le \hat u\}$ is nonempty.

In case of $\GG = \GGiso$ one can easily show that
$$
	\hat{\hat\ell}(x) \ = \ \sup_{t \in [0,x]} \, \hat\ell(t)	
 \quad\mbox{and}\quad	
	\hat{\hat u}(x)   \ = \ \inf_{s \in [x,1]} \, \hat u(s) .
$$
Note also that $\hat{\hat\ell}$ and $\hat{\hat u}$ are isotonic, whereas the raw functions $\hat\ell$ and $\hat u$ need not be.

\begin{figure}[h]
\centering
\includegraphics[height=7cm,width=10cm]{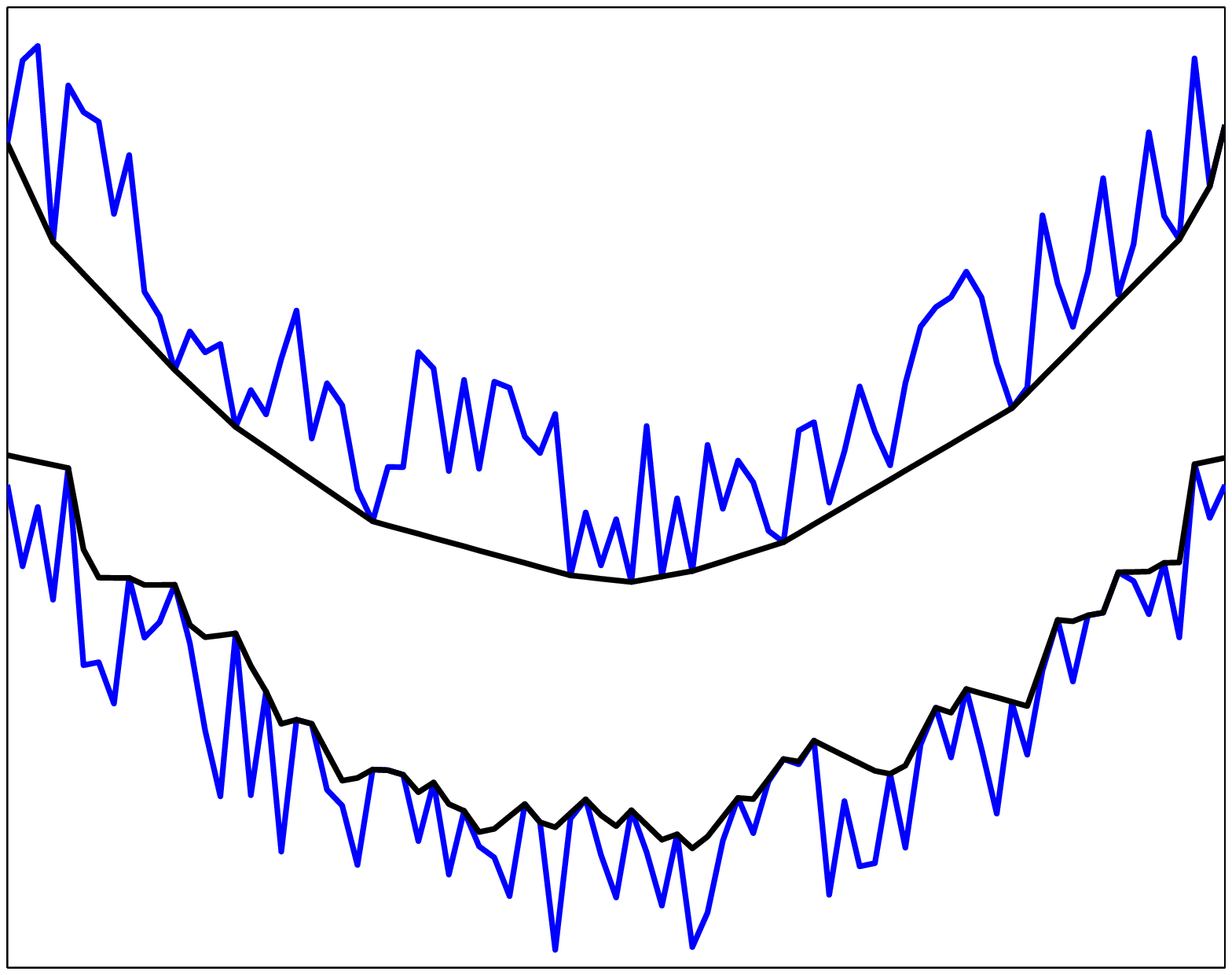}
\caption{Improvement $(\hat{\hat\ell},\hat{\hat u})$ of a band $(\hat\ell,\hat u)$ if $\GG = \GGconv$.}
\label{UpdateBand}
\end{figure}

In case of $\GG = \GGconv$ the modified upper bound $\hat{\hat u}$ is the greatest convex minorant of $\hat u$ and can be computed (in discrete models) by means of the pool-adjacent-violators algorithm (cf. Robertson et al.~1988). The modified lower bound $\hat{\hat\ell}(x)$ can be shown to be
$$
	\hat{\hat\ell}(x) 
	\ = \ \max \left\{ \sup_{0 \le s < t \le x} 
			\Bl( \hat{\hat u}(s) + {\hat\ell(t) - \hat{\hat u}(s)\over t-s} \, (x-s) \Br) , 
		\sup_{x \le s < t \le 1} 
			\Bl( \hat{\hat u}(t) - {\hat{\hat u}(t) - \hat\ell(s)\over t-s} \, (t-x) \Br) \right\} .
$$
This improved bound $\hat{\hat\ell}$ is not a convex function, though more regular than the raw function $\hat\ell$. Figure~\ref{UpdateBand} depicts some hypothetical confidence band $(\hat\ell,\hat u)$ for a function $f \in \GGconv$ and its improvement $(\hat{\hat\ell},\hat{\hat u})$.

%========================================
\subsection{Adaptivity in terms of rates}
%========================================

Whenever we construct a band following the recipe above we end up with a confidence band adapting to the unknown smoothness of $f$ in terms of rates of convergence. For $\beta, L > 0$ the H\"older smoothness class $\HH_{\beta,L}$ is defined as follows: In case of $0 < \beta \le 1$ let
$$
	\HH_{\beta,L} 
	\ := \ \Bl\{ g : |g(x) - g(y)| \le L |x-y|^\beta \mbox{ for all } x,y \Br\} .
$$
In case of $1 < \beta \le 2$ let
$$
	\HH_{\beta,L} \ := \ \Bl\{ g \in \CC^1 : g' \in \HH_{\beta-1,L} \Br\} .
$$

\begin{Theorem}	\label{Optimal Rates}
Suppose that $f \in \GG \cap \HH_{\beta,L}$, where either $\GG = \GGiso$ and $\beta \le 1$, or $\GG = \GGconv$ and $1 \le \beta \le 2$. Let $(\hat\ell,\hat u)$ be the confidence band for $f$ based on test functions $\psi^{(\ell)}, \psi^{(u)}$ as described previously. Then there exists a constant $\Delta$ depending only on $(\beta,L)$ and $(\psi^{(\ell)},\psi^{(u)})$ such that
$$
	\sup_{t \in [\epsilon_n,1-\epsilon_n]} \Bl( \hat u(t) - \hat\ell(t) \Br) 
	\ \le \ \Delta \rho_n \, 
		\Bl( 1 + {\kappa_\alpha + T(\psi^{(u)}) + T(-\psi^{(\ell)}) 
			\over \log(en)^{1/2}} \Br) ,
$$
where $\epsilon_n := \rho_n^{1/\beta}$ and
$$
	\rho_n \ := \ \Bl( {\log(en) \over n} \Br)_{}^{\beta/(2\beta+1)} .
$$
\end{Theorem}

Using the same arguments as Khas'minskii~(1978) one can show that for any $0\le r < s \le 1$,
$$
	\inf_{f \in \GG \cap \HH_{\beta,L}} \, 
	\Pr_f \left\{ \sup_{t \in [r,s]} (\hat u(t) - \hat\ell(t)) \le \Delta \rho_n \right\} 
	\ \to \ 0 ,
$$
provided that $\Delta > 0$ is sufficiently small. Thus our confidence bands adapt to the unknown smoothness of $f$.

%===========================================
\subsection{Testing hypotheses about $f(t)$}
%===========================================

In order to find suitable kernel functions $\psi^{(\ell)}, \psi^{(u)}$ we proceed similarly as D\"umbgen and Spokoiny~(2001, Section~3.2). That means we consider temporarily tests of the null hypothesis
$$
	\FF_o \ := \ \Bl\{ f \in \GG \cap L^2[0,1] : f(t) \le r-\delta \Br\}
$$
versus the alternative hypothesis
$$
	\FF_A \ := \ \Bl\{ f \in \GG \cap \HH_{k,L} : f(t) \ge r \Br\} .
$$
Here $t \in [0,1]$, $r \in \R$ and $L,\delta > 0$ are arbitrary fixed numbers, while
\be
	(\GG,k) \ =\ (\GGiso,1)	
	\quad\mbox{or}\quad	
	(\GG,k) \ = \ (\GGconv,2) .
	\label{GGk}
\ee
Note that $\FF_o$ and $\FF_A$ are closed, convex subsets of $L^2[0,1]$. Suppose that there are functions $f_o\in\FF_o$ and $f_A\in\FF_A$ such that
$$
	\int_0^1 (f_o-f_A)(x)^2 \, dx 
	\ = \ \min_{g_o\in\FF_o, \, g_A\in\FF_A} \int_0^1 (g_o-g_A)(x)^2 \, dx .
$$
Then optimal tests of $\FF_o$ versus $\FF_A$ are based on the linear test statistic
$\int_0^1 (f_A-f_o) \, dY$, where critical values have to be computed under the assumption $f=f_o$. The problem of finding such functions $f_o,f_A$ is treated in Section~\ref{Optimization}. Here is the conclusion: Let
\be
	\psi^{(\ell)}(x) \ := \ \left\{\barr{ll} 
		1\{x \in [-1,0]\} (1 + x) 
			& \mbox{if } \GG = \GGiso , \\
		1\{x \in [-2,2]\} \Bl( 1 - (3/2)|x| + x^2/2 \Br) 
			& \mbox{if } \GG = \GGconv .
	\earr\right.
	\label{PsiL}
\ee
Then the functions
\be
	f_A(s) \ := \ \left\{\barr{cl}
		r + L(s-t)     & \mbox{if } \GG = \GGiso \\
		r + L(s-t)^2/2 & \mbox{if } \GG = \GGconv 
	\earr\right.
	\label{Least favorable fA}
\ee 
and 
$$
	f_o \ := \ f_A - \delta \psi^{(\ell)}_{h,t}	
	\quad\mbox{with } h := (\delta/L)^{1/k}
$$
solve our minimzation problem, provided that $a^{(\ell)}h \le t \le 1 - b^{(\ell)}h$. Thus the optimal linear test statistic may be written as $\int_0^1 \psi_{h,t} \, dY = \psi Y(h,t)$. Elementary considerations show that the inequality
$$
	\hat f_h^{(\ell)}(t) 
		- {\|\psi^{(\ell)}\| (\Gamma(d^{(\ell)} h) + \kappa_\alpha)
			\over \langle 1,\psi^{(\ell)}\rangle (nh)^{1/2}} 
	\ \le \ r_o
$$
is equivalent to
\bea
	\psi Y(h,t) 
	& \le & n^{1/2} h r_o \langle 1, \psi^{(\ell)}\rangle 
		+ h^{1/2} \|\psi^{(\ell)}\| (\Gamma(d^{(\ell)} h) + \kappa_\alpha) \\
	& = & \Ex_{f_o}(\psi Y(h,t)) 
		+ \Var(\psi Y(h,t))^{1/2} (\Gamma(d^{(\ell)} h) + \kappa_\alpha) .
\eea
Thus our lower confidence bound $\hat\ell$ may be interpreted as a multiple test of all null hypotheses $\{f \in \GG : f(t) \le r_o\}$ with $t \in [0,1]$ and $r_o \in \R$.

Analogous considerations yield a candidate for $\psi^{(u)}$: Let
$$
	\FF_o \ := \ \Bl\{ f \in \GG \cap L^2[0,1] : f(t) \ge r+\delta \Br\}
$$
and
$$
	\FF_A \ := \ \Bl\{ f \in \GG \cap \HH_{k,L} : f(t) \le r \Br\} .
$$
Then the function $f_A$ in (\ref{Least favorable fA}) and
$$
	f_o \ := \ f_A + \delta \psi^{(u)}_{h,t}	
	\quad\mbox{with } h := (\delta/L)^{1/k}
$$
form a least favorable pair $(f_o,f_A)$ in $\FF_o \times \FF_A$, where
\be
	\psi^{(u)}(x) \ := \ \left\{\barr{ll} 
		1\{x \in [0,1]\} (1 - x) 
			& \mbox{if } \GG = \GGiso , \\
		1\{x \in [-2^{1/2},2^{1/2}]\} (1 - x^2/2) 
			& \mbox{if } \GG = \GGconv .
	\earr\right.
	\label{PsiU}
\ee

\begin{figure}[h]
\includegraphics[height=6cm,width=7cm]{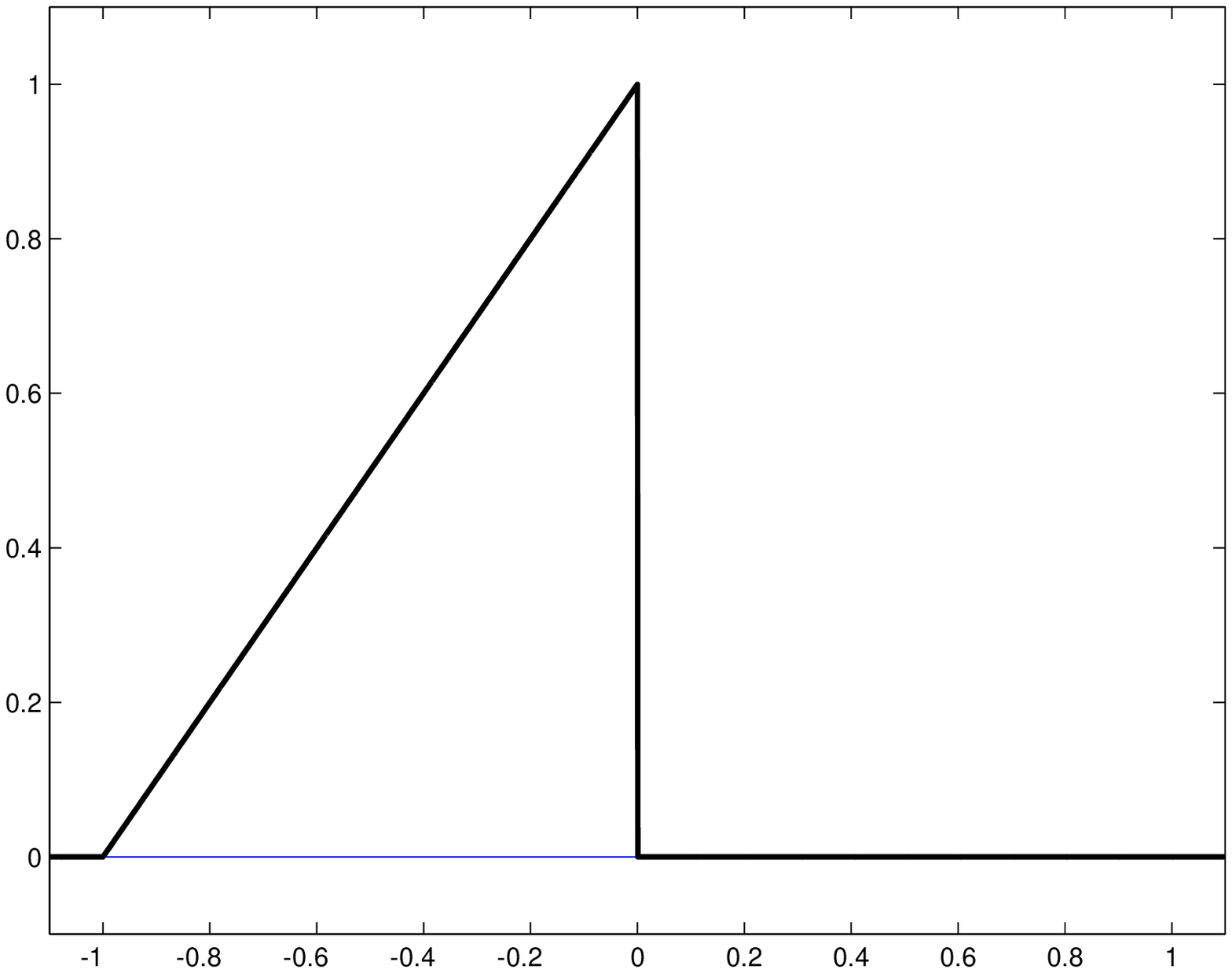}
\hfill 
\includegraphics[height=6cm,width=7cm]{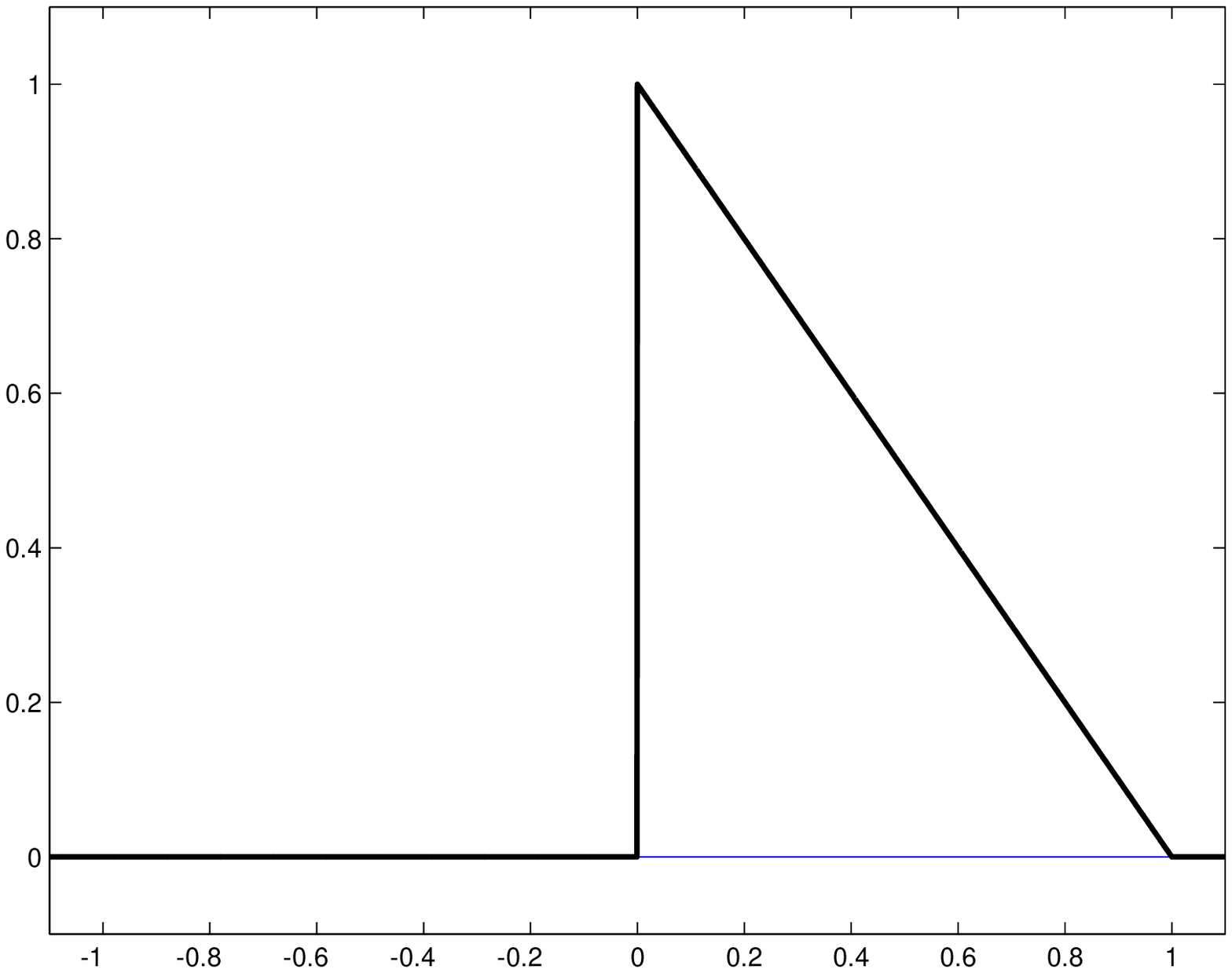}
\caption{Kernel functions $\psi^{(\ell)}, \psi^{(u)}$ for $\GGiso$.}
\label{OptKernelsGGiso}
\end{figure}

\begin{figure}[h]
\includegraphics[height=6cm,width=7cm]{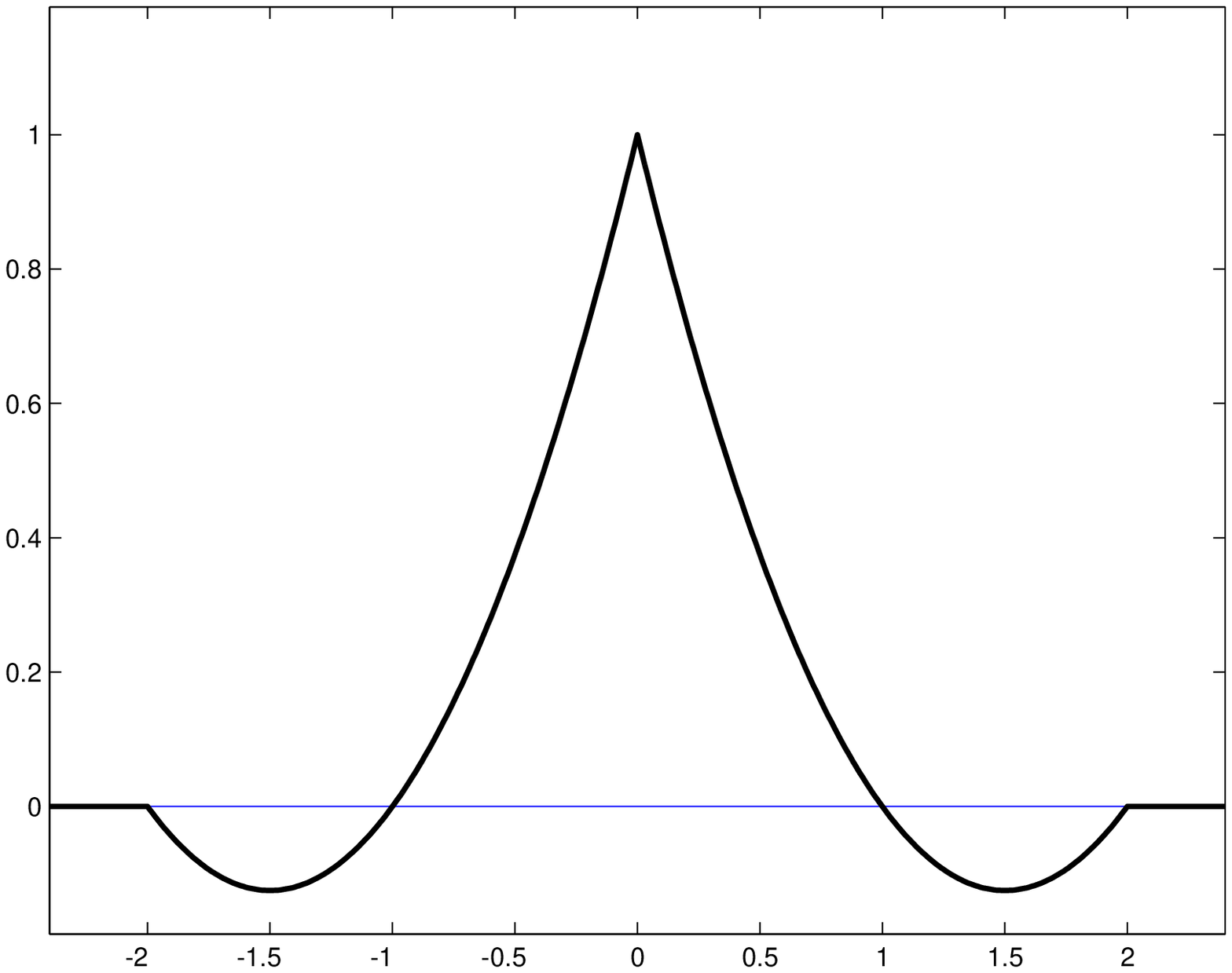}
\hfill 
\includegraphics[height=6cm,width=7cm]{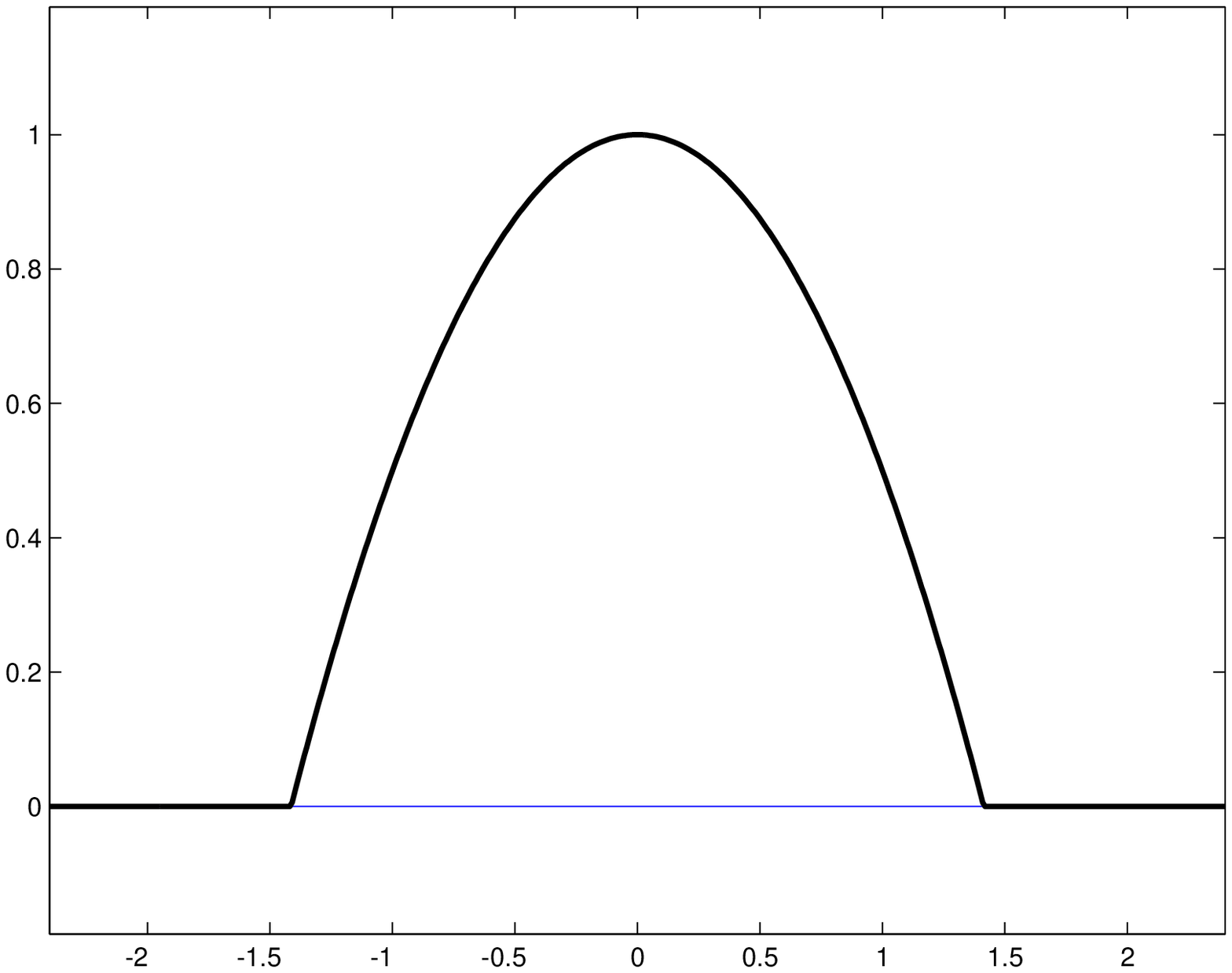}
\caption{Kernel functions $\psi^{(\ell)}, \psi^{(u)}$ for $\GGconv$.}
\label{OptKernelsGGconv}
\end{figure}

Figures~\ref{OptKernelsGGiso} and \ref{OptKernelsGGconv} depict the functions $\psi^{(\ell)}$ in (\ref{PsiL}) and $\psi^{(u)}$ in (\ref{PsiU}).

%======================================================
\subsection{Optimal constants and local adaptivity}
%======================================================

Now we are going to show that our multiscale confidence band $(\hat\ell, \hat u)$, if constructed with the kernel functions in (\ref{PsiL}) and (\ref{PsiU}), is locally adaptive in a certain sense. Precisely, we consider an arbitrary fixed function $f_o \in \GG\cap\CC^k$ with $(\GG,k)$ as specified in (\ref{GGk}). We analyze quantities such as
$$
	\|(\hat u - f_o) w\|^{+}_{r,s}	
	\quad\mbox{and}\quad	
	\|(f_o - \hat\ell) w\|^{+}_{r,s} ,
$$
where $w$ is some positive weight function on the unit interval and
$$
	\|g\|^{+}_{r,s} \ := \ \sup_{t \in [r,s]} \, g(t) .
$$
The function $w$ should reflect local smoothness properties of $f_o$ in an appropriate way. The following theorem demonstrates that the $k$--th derivative of $f_o$, denoted by $\Dk f_o$, plays a crucial role.

\begin{Theorem}	\label{Adapt I}
For arbitrary fixed numbers $0 \le r < s \le 1$ let
$$
	L \ := \ \max_{t \in [r,s]} \, \Dk f_o(t) .
$$
Then for any $\gamma\in\left]0,1\right[$,
\bea
	\inf_{(\hat\ell,\hat u)} \, \Pr_{f_o} \Bl\{ \|f - \hat\ell\|^{+}_{r,s} 
		\ge \gamma \Delta^{(\ell)} L^{1/(2k+1)} \rho_n \Br\} 
	& \ge & 1 - \alpha + o(1) , \\
	\inf_{(\hat\ell,\hat u)} \, \Pr_{f_o} \Bl\{ \|\hat u - f\|^{+}_{r,s} 
		\ge \gamma \Delta^{(u)}    L^{1/(2k+1)} \rho_n \Br\} 
	& \ge & 1 - \alpha + o(1) ,
\eea
where both infima are taken over all confidence bands $(\hat\ell,\hat u)$ satisfying (\ref{Confidence}), and
\bea
	\Delta^{(z)} & := & \Bl( (k+1/2) \|\psi^{(z)}\|^2 \Br)^{-k/(2k+1)} , \\
	\rho_n     & := & \Bl( {\log(en)\over n} \Br)^{k/(2k+1)} .
\eea
\end{Theorem} 

In case of $\GG=\GGiso$, the critical constants are $\Delta^{(\ell)} = \Delta^{(u)} = 2^{1/3} \approx 1.260$. In case of $\GG=\GGconv$,
$$
	\Delta^{(\ell)} \ = \ (3/4)^{2/5} \ \approx \ 0.891	
	\quad\mbox{and}\quad	
	\Delta^{(u)} \ = \ 3^{2/5}/128^{1/5} \ \approx \ 0.588 .
$$
This indicates that bounding a convex function from below is more difficult than finding an upper bound.

In view of Theorem~\ref{Adapt I} we introduce for arbitrary fixed $\epsilon > 0$ the weight function
$$
	w_\epsilon \ := \ \Bl( \max(\Dk f_o,\epsilon) \Br)^{-1/(2k+1)}
$$
reflecting the local smoothness of $f_o$. The next theorem shows that our particular confidence band $(\hat\ell,\hat u)$ attains the lower bounds of Theorem~\ref{Adapt I} pointwise. Suprema such as $\|(f_o - \hat\ell) w_\epsilon\|^{+}_{r,s}$ and $\|(\hat u - f_o) w_\epsilon\|^{+}_{r,s}$ attain their respective lower bounds $\Delta^{(\ell)}$, $\Delta^{(u)}$ up to a multiplicative factor $2^{k/(k+1/2)} + o_p(1)$.

\begin{Theorem}	\label{Adapt II}
Let $(\hat\ell,\hat u)$ be the confidence band based on the kernel functions in (\ref{PsiL}) and (\ref{PsiU}). If $f = f_o$, then for arbitrary $\epsilon > 0$ and any $t\in\left]0,1\right[$,
\bea
	(f_o - \hat\ell)(t) w_\epsilon(t) 
	& \le & \left( \Delta^{(\ell)} + o_p(1) \right) \rho_n , \\
	(\hat u - f_o)(t) w_\epsilon(t) 
	& \le & \left( \Delta^{(u)}    + o_p(1) \right) \rho_n .
\eea
Moreover,
\bea
	\|(f_o - \hat\ell) w_\epsilon\|^{+}_{\epsilon,1-\epsilon} 
	& \le & \left( 2^{k/(k+1/2)} \Delta^{(\ell)} + o_p(1) \right) \rho_n , \\
	\|(\hat u - f_o) w_\epsilon \|^{+}_{\epsilon,1-\epsilon} 
	& \le & \left( 2^{k/(k+1/2)} \Delta^{(u)} + o_p(1) \right) \rho_n .
\eea
\end{Theorem} 

If we used kernel functions differing from (\ref{PsiL}) and (\ref{PsiU}), then pointwise optimality would be lost, and the constants for the supremum distances would get worse.

%===========================================
\section{Simulations and numerical examples}
\label{Examples}
%===========================================

Here we demonstrate the performance of the procedures in Section~\ref{GG smooth}. We replace the continuous white noise model with a discrete one: Suppose that one observes a random vector $\vec{Y} \in \R^n$ with components
\be
	Y_i \ = \ f(x_i) + \epsilon_i ,
	\label{WN Regression}
\ee
where $x_i := (i-1/2)/n$, and the random errors $\epsilon_i$ are independent
with Gaussian distribution $\NN(0,\sigma^2)$. Our kernel functions $\psi^{(\ell)}$ and $\psi^{(u)}$ are rescaled as follows:
\bea
	\psi^{(\ell)}(x) 
	& := & \left\{\barr{ll} 
		1\{x \in [-1,0]\} (1 + x) 
			& \mbox{if } \GG = \GGiso , \\
		1\{x \in [-1,1]\} \Bl( 1 - 3 |x| + 2 x^2 \Br) 
			& \mbox{if } \GG = \GGconv ,
	\earr\right. \\
	\psi^{(u)}(x) 
	& := & \left\{\barr{ll} 
		1\{x \in [0,1]\} (1 - x) 
			& \mbox{if } \GG = \GGiso , \\
		1\{x \in [-1,1]\} (1 - x^2) 
			& \mbox{if } \GG = \GGconv .
	\earr\right.
\eea
Note that now $a^{(\ell)}, a^{(u)}, b^{(\ell)}, b^{(u)} \in \{0,1\}$. For convenience we compute kernel estimators and confidence bounds for $f$ only on the grid $\TT_n:=\{1/n,2/n,\ldots,1 - 1/n\}$, while the bandwidth parameter $h$ is restricted to
$$
	H_n \ := \ \left\{\barr{cl} 
		\{1/n, 2/n, \ldots, 1\} & \mbox{if } \GG = \GGiso , \\
		\{1/n, 2/n, \ldots, \lfloor n/2\rfloor/n\} & \mbox{if } \GG = \GGconv .
	\earr\right.
$$
Let $\psi$ stand for $\psi^{(\ell)}$ or $\psi^{(u)}$ with support $[-a,b]$. Then for $h\in H_n$ and $t\in \TT_n$ with $ah \le t \le 1 - bh$ we define
$$
	\psi\vec{Y}(h,t) 
	\ := \ \sum_{i=1}^n \psi \Bl( {x_i - t\over h} \Br) Y_i 
	\ =  \ \sum_{j=1 - anh}^{bnh} 
		\psi \Bl( {j - 1/2\over nh} \Br) Y_{nt+j}
$$
and
$$
	\hat f_h(t) \ := \ {\psi\vec{Y}(h,t) \over S_{nh}} ,
$$
where $S_d$ stands for $\sum_{j=1-d}^d \psi((j-1/2)/d)$. The standard deviation of $\hat{f}_h(t)$ equals $\sigma_h := \sigma R_{nh}^{1/2}/S_{nh}$, where $R_d := \sum_{j=1-d}^d \psi((j-1/2)/d)^2$. Tedious but elementary calculations show that in case of $\GG = \GGiso$,
$$
	S_d \ = \ d/2	\quad\mbox{and}\quad	R_d \ = \ d/3 - 1/(12d) .
$$
In case of $\GG = \GGconv$,
$$
	\barr{ccccccc}
	S^{(\ell)}_d & = & d/3 - 1/(3d) & \mbox{and} & 
	R^{(\ell)}_d & = & 4d/15 - 1/(2d) + 7/(30d^3) , \\
	S^{(u)}_d    & = & 4d/3 + 1/(6d) & \mbox{and} & 
	R^{(u)}_d    & = & 16d/15 + 7/(120d^3) .
	\earr
$$
Note that here $S^{(\ell)}_1 = 0 = \psi^{(\ell)}\vec{Y}(1/n,\cdot)$, whence the bandwidth $1/n$ is excluded from any computation involving $\psi^{(\ell)}$.

As for the bias of these kernel estimators, one can deduce from Lemma~\ref{Ortho} that
$\Ex \hat f^{(\ell)}_h(t) \le f(t)$ and $\Ex \hat f^{(u)}_h(t) \ge f(t)$ whenever $f\in\GG$. Here is a discrete version of our multiscale test statistic: $T_n^* := \max \Bl( T_n(\psi^{(\ell)}), T_n(- \psi^{(u)}) \Br)$, where
$$
	T_n(\pm\psi) \ := \ \max_{h\in H_n} \ \max_{t\in \TT_n \cap [ah,1-bh]} 
		\Bl( \pm \sigma^{-1} R_{nh}^{-1/2} \psi\vec{E}(h,t) - \Gamma((a+b)h) \Br)
$$
with $\vec{E} := (\epsilon_i)_{i=1}^n$. Let $\kappa_{\alpha,n}$ be the $(1-\alpha)$--quantile of $T_n^*$. Then
\bea
	\hat\ell(t) 
	& := & \max_{h\in H_n \, : \, t \in [a^{(\ell)}h, 1 - b^{(\ell)}h]} 
		\Bl( \hat{f}^{(\ell)}_h(t)
			- \sigma^{(\ell)}_h (\Gamma(d^{(\ell)} h) + \kappa_{\alpha,n}) \Br) , \\
	\hat{u}(t) 
	& := & \min_{h\in H_n \, : \, t \in [a^{(u)}h, 1 - b^{(u)}h]} 
		\Bl( \hat{f}^{(u)}_h(t)
			+ \sigma^{(u)}_h (\Gamma(d^{(u)} h) + \kappa_{\alpha,n}) \Br) ,
\eea
defines a confidence band for $f$ such that
$$
	\Pr \Bl\{ \hat\ell \le f \le \hat{u} \mbox{ on } \TT_n \Br\} \ \ge \ 1 - \alpha
	\quad\mbox{whenever } f \in \GG .
$$
Equality holds if $\GG = \GGiso$ and $f$ is constant, or if $\GG = \GGconv$ and $f$ is linear. If the noise variance $\sigma^2$ is unknown, it may be estimated as described in D\"umbgen and Spokoiny~(2001). Then, under moderate regularity assumptions on $f$, our confidence bands have {\sl asymptotic} coverage probability at least $1 - \alpha$ as $n$ tends to infinity.

{\bf Critical values.} For various values of $n$ we estimated several quantiles $\kappa_{\alpha,n}$ in 9999 Monte-Carlo simulations; see Table~\ref{Simulations}. One can easily show that the critical value $\kappa_{\alpha,n}$ converges to the corresponding quantile $\kappa_\alpha$ for the continuous white noise model as $n \to \infty$. Software for the computation of critical values as well as confidence bands may be obtained from the author's URL.

\begin{table}[h]
\centerline{\begin{tabular}{|c||c|c|c||c|c|c|}	\hline
	     & \multicolumn{3}{|c||}{$\GGiso$}
	     & \multicolumn{3}{c|}{$\GGconv$} \\
	 $n$ & $\kappa_{0.5,n}$ & $\kappa_{0.1,n}$ & $\kappa_{0.05,n}$
	     & $\kappa_{0.5,n}$ & $\kappa_{0.1,n}$ & $\kappa_{0.05,n}$ \\\hline\hline
	 100 & 0.330 & 1.092 & 1.349 & 0.350 & 1.053 & 1.283 \\\hline
	 200 & 0.433 & 1.146 & 1.392 & 0.430 & 1.121 & 1.342 \\\hline
	 300 & 0.475 & 1.169 & 1.416 & 0.470 & 1.126 & 1.342 \\\hline
	 400 & 0.507 & 1.204 & 1.446 & 0.489 & 1.128 & 1.340 \\\hline
	 500 & 0.526 & 1.222 & 1.450 & 0.512 & 1.143 & 1.358 \\\hline
	 700 & 0.570 & 1.252 & 1.492 & 0.536 & 1.162 & 1.380 \\\hline
	1000 & 0.585 & 1.250 & 1.483 & 0.552 & 1.178 & 1.393 \\\hline
\end{tabular}}
\caption{Some critical values for the discrete white noise model}
\label{Simulations}
\end{table}

\begin{figure}
\centering
\includegraphics[height=9cm,width=12cm]{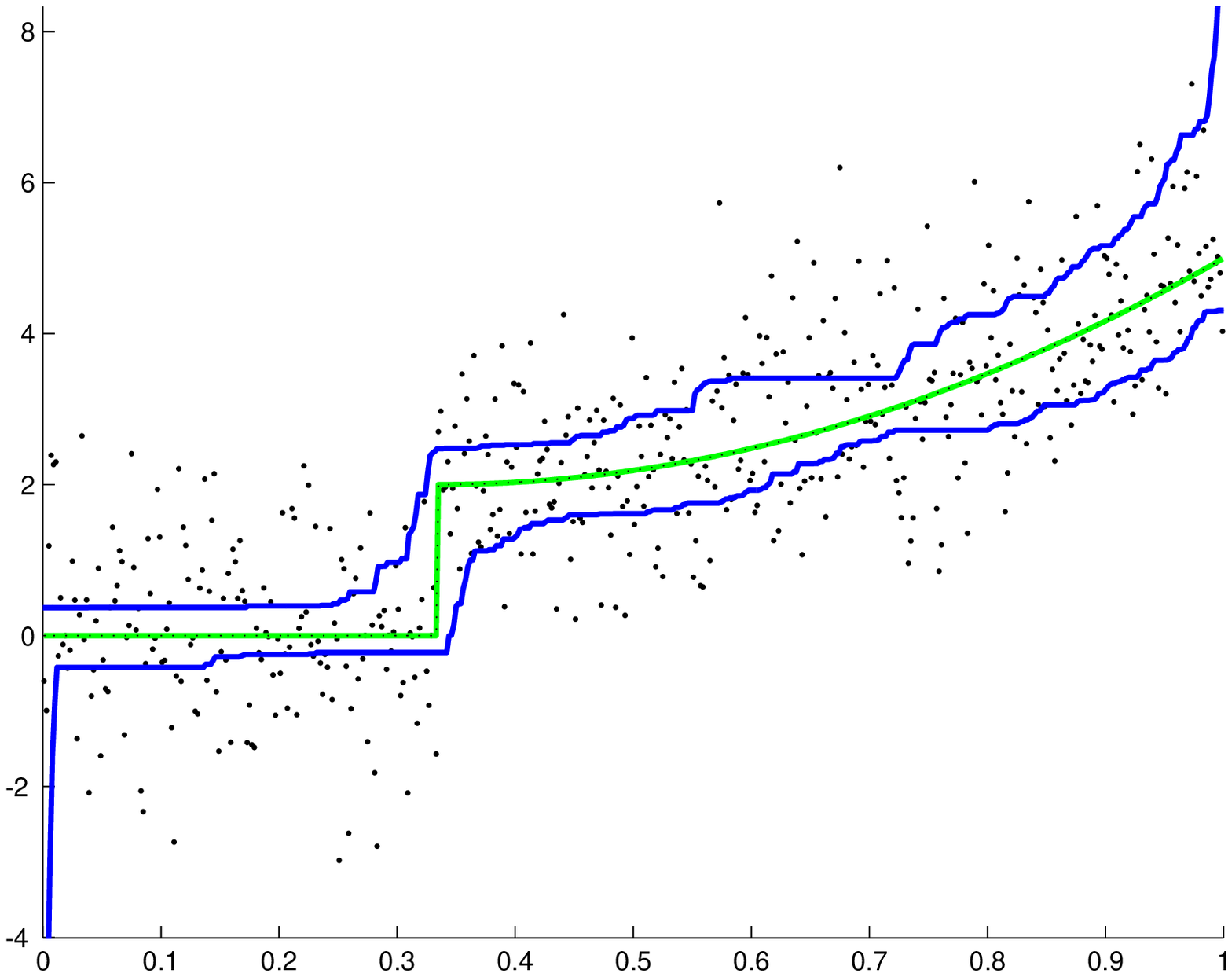}
\caption{Data $\vec{Y}$ and $95\%$--confidence band for $f \in \GGiso$.}
\label{ExampleIso}
\end{figure}

{\bf Two numerical examples.} Figure~\ref{ExampleIso} shows a simulated data vector $\vec{Y}$ with $n = 500$ components together with the corresponding $95\%$--confidence band $(\hat\ell,\hat u)$ after postprocessing, where $f$ is assumed to be isotonic. The latter function is depicted as well. Note that the band is comparatively narrow in the middle of $\left]0,1/3\right[$, on which $f$ is constant. On $\left]1/3,1\right]$ the width $\hat u - \hat\ell$ tends to inrease, as does $\nabla f$. These findings are in accordance with Theorem~\ref{Adapt II}.

\begin{figure}
\centering
\includegraphics[height=9cm,width=12cm]{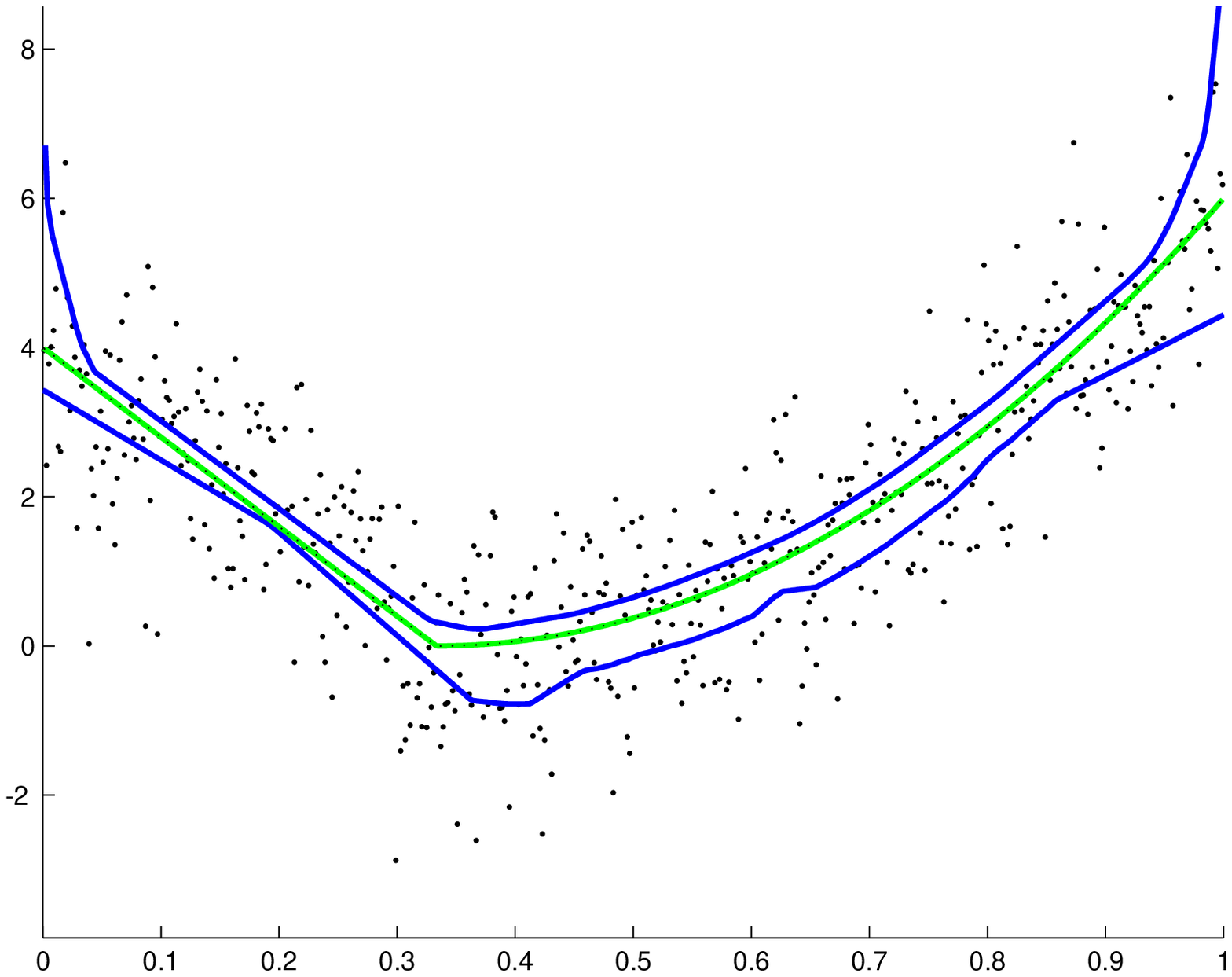}
\caption{Data $\vec{Y}$ and $95\%$--confidence band for $f \in \GGconv$.}
\label{ExampleConv}
\end{figure}

An analogous plot for a convex function $f$ can be seen in Figure~\ref{ExampleConv}. Note that the deviation $f - \hat\ell$ is mostly greater than $\hat u - f$, as predicted by Theorem~\ref{Adapt II}.

%================
\section{Proofs}
\label{Proofs}
%================

\noindent
{\bf Proof of Theorem~\ref{Lower Bounds GGiso}.} In order to prove lower bounds we 
construct unfavorable subfamilies of $\GGiso$ similarly as Khasminski~(1978). For a given integer $m > 0$ we define $I_1 := [0,1/m]$ and $I_j := \left](j-1)/m, j/m\right]$ for $1 < j \le m$. Then we define step functions $g$ and $h^{}_\xi$ for $\xi \in \R^m$ via
$$
	g(t) \ := \ 2j - 1	\quad\mbox{and}\quad	
	h^{}_\xi(t) \ := \ \xi_j	\quad\mbox{for } t \in I_j, 1 \le j \le m .
$$
For any $\delta > 0$ and $\xi \in [-\delta, \delta]^m$ the function $\delta g + h^{}_\xi$ is isotonic on $[0,1]$. Now we restrict our attention to the parametric submodel $\FF_o = \Bl\{ \delta g + h^{}_\xi : \xi \in [-\delta,\delta]^m \Br\}$ of $\GGiso \cap L^2[0,1]$. Any confidence band $(\hat\ell,\hat u)$ for $f = \delta g + h^{}_\xi$ defines a confidence set $S = S_1 \times S_2 \times \cdots \times S_m$ for $\xi$ via
$$
	S_j 
	\ := \ \Bl[ \sup_{t\in I_j} \, \hat\ell(t) - \delta (2j-1) , 
	            \inf_{t\in I_j} \, \hat   u(t) - \delta (2j-1) \Br] .
$$
Here $\hat\ell \le f \le \hat u$ if, and only if, $\xi \in S$. Moreover,
$$
	D_\epsilon(\hat\ell, \hat u) \ \ge \ \max_{j=1,\ldots,m} {\rm length}(S_j)	
	\quad\mbox{for } 1/(m+1) \le \epsilon < 1/m .
$$
However,
\bea
	\log {d\Pr_{\delta g + h^{}_\xi} \over d\Pr_{\delta g}} (Y) 
	& = & n^{1/2} \int_0^1 h^{}_\xi \, d\tilde{Y} - n \int_0^1 h^{}_\xi(t)^2 \, dt / 2 \\
	& = & \sum_{j=1}^m \Bl( (n/m)^{1/2} \xi_j X_j - (n/m) \xi_j^2/2 \Br) \\
	& = & \log {d\NN((n/m)^{1/2}\xi,I)\over d\NN(0,I)} (X) ,
\eea
where $\tilde{Y}(t) := Y(t) - n^{1/2} \int_0^t \delta g(s) \, ds$ and $X := (X_j)_{j=1}^m$ with components
$$
	X_j 
	\ := \ m^{1/2} \Bl( \tilde{Y}(j/m) - \tilde{Y}((j-1)/m) \Br) .
$$
In case of $f = \delta g$ these random variables are independent and standard normal. Consequently, $X$ is a sufficient statistic for the parametric submodel $\FF_o$ with distribution $\NN_m((n/m)^{1/2}\xi,I)$ in case of $f = \delta g + h^{}_\xi$. In particular, the conditional distribution of $S$ given $X$ does not depend on $\xi$. Hence letting $\delta = (n/m)^{-1/2} c_m$ with $c_m := (2\log m)^{1/2}$ it follows from Theorem~\ref{CS Sup-Diameter}~(b) in Section~\ref{Decision Theory} that for $1/(m+1) \le \epsilon < 1/m$,
\bea
	\lefteqn{ \inf_{f \in \GGiso \cap L^2[0,1]} \, 
		\Pr_f \Bl\{ \hat\ell \le f \le \hat u \mbox{ and } D_\epsilon(\hat\ell,\hat u) 
			\le 2 {c_m - b_m \over (n/m)^{1/2}} \Br\} } \\
	& \le & \min_{\xi \in [-\delta, \delta]^m} \, 
		\Pr_\xi \Bl\{ \xi \in S \mbox{ and } \max_{j=1,\ldots,m} {\rm length}(S_j) 
			\le 2 {c_m - b_m \over (n/m)^{1/2}} \Br\} 
		\ \le \ b_m ,
\eea
where $b_1, b_2, b_3, \ldots$ are universal positive numbers such that $\lim_{m\to\infty} b_m = 0$. This entails the assertion of Theorem~\ref{Lower Bounds GGiso} with $\log(1/\epsilon)$ in place of $\log(e/\epsilon)$ and
$$
	b(\epsilon) \ := \ (2 \log(1/\epsilon))^{1/2} - (m\epsilon)^{1/2} (c_m - b_m)	
	\quad\mbox{for } 1/(m+1) \le \epsilon < 1/m .
$$
Finally note that $\log(e/\epsilon)^{1/2} = \log(1/\epsilon)^{1/2} + o(1)$ as $\epsilon \downarrow 0$.	\hfill	$\Box$

\noindent
{\bf Proof of Theorem~\ref{Optimal Rates}.} Instead of an upper bound for $\hat u-\hat\ell$ we prove an upper bound for $\hat u-f$, because analogous arguments apply to $f-\hat\ell$. In what follows let $\psi = \psi^{(u)}$ with support $[-a,b]$. For $t \in [0,1]$ and $h>0$ with $ah\le t\le 1-bh$,
\bean
	\hat u(t) - f(t) 
	& \le & \hat f_h(t) - f(t) 
		+ {\|\psi\| (\Gamma((a+b)h)+\kappa_\alpha) \over \langle 1,\psi\rangle (nh)^{1/2}} 
	\nonumber \\
	& =    & {\Bl\langle f(t + h\,\cdot) - f(t),\psi \Br\rangle \over \langle 1,\psi\rangle} 
		+ {\psi W(h,t) \over n^{1/2} h \langle 1,\psi\rangle} 
		+ {\|\psi\| (\Gamma((a+b)h)+\kappa_\alpha) \over \langle 1,\psi\rangle (nh)^{1/2}} 
	\nonumber \\
	& \le & {\Bl\langle f(t + h\,\cdot) - f(t),\psi \Br\rangle \over \langle 1,\psi\rangle} 
		+ {\|\psi\| \Bl( 2 \Gamma((a+b)h) + \kappa_\alpha + T(\psi) \Br) 
			\over \langle 1,\psi\rangle (nh)^{1/2}} .
	\label{Optimal Rates.1}
\eean
For any function $g \in \HH_{\beta,L}$,
\bea
	\Bl| g(x) - g(0) \Br| & \le    & L |x|^\beta	
		\quad\mbox{if } \beta \le 1 , \\
	\Bl| g(x) - g(0) - g'(0) x \Br| & \le & L |x|^\beta	
		\quad\mbox{if } 1 < \beta \le 2 .
\eea
Since $f(t + h\,\cdot) \in \HH_{\beta,L h^\beta}$ if $f \in \HH_{\beta,L}$, this implies that
$$
	{\Bl\langle f(t + h\,\cdot) - f(t),\psi \Br\rangle \over \langle 1,\psi\rangle} 
	\ \le \ {L h^\beta \int_{-a}^b |x|^\beta |\psi(x)| \, dx \over \langle 1,\psi\rangle} 
	\ \le \ \Delta h^\beta .
$$
Here and subsequently $\Delta$ denotes a generic constant depending only on $(\beta,L)$ and $\psi$. Its value may vary from one place to another. In case of $t \in [\epsilon_n,1-\epsilon_n]$ and $h = \epsilon_n/\max(a,b)$ the right-hand side of (\ref{Optimal Rates.1}) is not greater than
$$
	\Delta \epsilon_n^\beta 
	+ {\Delta \Bl( \log(en)^{1/2} + \kappa_\alpha + T(\psi) \Br) \over (n\epsilon_n)^{1/2}} 
	\ = \ \Delta \rho_n \Bl( 1 + {\kappa_\alpha + T(\psi)\over \log(en)^{1/2}} \Br) .
	\eqno{\Box}
$$

\noindent
{\bf Proof of Theorem~\ref{Adapt I}.} We prove only the lower bound for $f_o - \hat\ell$, because $\hat u - f_o$ can be treated analogously. It suffices to consider the case $L > 0$ and to show that for any fixed number $\gamma \in \left]0,1\right[$,
$$
	\Pr_{f_o} \Bl\{ \|f_o - \hat\ell\|^{+}_{r,s} 
		\ge \gamma \Delta^{(\ell)} L^{1/(2k+1)} \rho_n \Br\} 
	\ \ge \ 1 - \alpha + o(1)
$$
for arbitrary confidence bands $(\hat\ell,\hat u) = (\hat\ell_n,\hat u_n)$ satisfying (\ref{Confidence}). Without loss of generality one may assume that
$$
	\Dk f_o \ \ge \ L	\quad\mbox{on } [r,s] .
$$
Otherwise one could increase $\gamma$ and decrease $L$ without changing $\gamma L^{1/(2k+1)}$, and replace $[r,s]$ with some nondegenerate subinterval. Let $\psi$ stand for $\psi^{(\ell)}$ with support $[-a,b]$. For $0 < h \le (s-r)/(a+b)$ and positive integers $j\le m:=\lfloor (s-r)/((a+b)h\rfloor$ let
$$
	t_j \ := \ s + ah + (j-1)(a+b) h	
	\quad\mbox{and}\quad	
	f_j \ := \ f_o - L h^k \psi^{}_{h,t_j} .
$$
It follows from Lemma~\ref{Nice Psis} that these functions $f_j$ belong to $\GG\cap L^2[0,1]$. Thus (\ref{Confidence}) implies that the event
$$
	A \ := \ \Bl\{ \hat\ell \le f_j \mbox{ for some } j \le m \Br\}
$$
satisfies the inequality $\Pr_{f_j}(A) \ge 1-\alpha$ for all $j \le m$. Since $\|f_o - f_j\|^{+}_{r,s} \ge \delta$, this entails the inequality
$$
	\Pr_{f_o} \Bl\{ \|f_o - \hat\ell\|^{+}_{r,s} \ge L h^k \Br\} 
	\ \ge \ \Pr_{f_o}(A) 
	\ \ge \ 1 - \alpha - \min_{j \le m} \Bl (\Pr_{f_j}(A) - \Pr_{f_o}(A) \Br) .
$$
Now let $h := (c\rho_n)^{1/k}$ so that $Lh^k=Lc\rho_n$, where $c > 0$ is some number to be specified later. For sufficiently large $n$ this bandwidth $h$ is smaller than $(s-r)/(a+b)$. Then
$$
	\log {d\Pr_{f_j}\over d\Pr_{f_o}}(Y) 
	\ = \ n^{1/2} h^{k+1/2} L \|\psi\| X_j - nh^{2k+1} L^2 \|\psi\|^2/2 ,
$$
where $X_j := h^{-1/2} \|\psi\|^{-1} \int_0^1 \psi_{h,t_j} \, d\tilde{Y}$ and $\tilde{Y}(t) := Y(t) - n^{1/2} \int_0^t f_o(x) \, dx$. Thus $X := (X_j)_{j=1}^m$ is a sufficient statistic for the restricted model $\{f_o, f_1, f_2, \ldots, f_m\}$, where $\LL_{f_o}(X)$ is a standard normal distribution on $\R^m$. Thus it follows from Theorem~\ref{CS Sup-Diameter}~(a) and a standard sufficiency argument that
$$
	\lim_{n\to\infty} \, \min_{1 \le j \le m} \Bl (\Pr_{f_j}(A) - \Pr_{f_o}(A) \Br) \ = \ 0	
	\quad\mbox{if}\quad	
	\lim_{n\to\infty} \, {nh^{2k+1} L^2 \|\psi\|^2 \over 2 \log m} \ < \ 1 .
$$
Since $\log m = (1 + o(1)) \log(n)/(2k+1)$, the limit on the right hand side is equal to
$$
	c^{(2k+1)/k} L^2 \|\psi\|^2 (k+1/2)
$$
and smaller than one if $c$ equals $\gamma\Delta^{(\ell)}L^{-2k/(2k+1)}$. In that case, the lower bound $Lh^k=Lc\rho_n$ for $\|f_o-\hat\ell\|^{+}_{r,s}$ equals $\gamma \Delta^{(\ell)} L^{1/(2k+1)} \rho_n$ as desired.	\hfill	$\Box$

\noindent
{\bf Proof of Theorem~\ref{Adapt II}.} Again we restrict our attention to $f_o - \hat\ell$ and let $\psi := \psi^{(\ell)}$ with support $[-a,b]$. For any fixed $\epsilon > 0$ and arbitrary $t \in [0,1]$ let $h_t > 0$ and
$$
	L_t \ := \ \max_{s \in [t-ah_t, t+bh_t] \cap [0,1]} \, \max(\Dk f_o(s),\epsilon) .
$$
In case of $ah_t \le t \le 1 - bh_t$ the inequality $(f_o - \hat\ell)(t) \ge L_t h_t^k$ implies that 
$$
	\hat f_{h_t}(t) 
		- { \|\psi\| \Bl( \Gamma((a+b)h_t) + \kappa_\alpha \Br) 
			\over (nh_t)^{1/2} \langle 1,\psi\rangle } 
	\ \le \ f_o(t) - L_t h_t^k .
$$
Since $f=f_o$, this can be rewritten as
\bea
	{\psi W(h_t,t) \over h_t^{1/2} \|\psi\|} 
	& \le & - \, {(nh_t)^{1/2} \over \|\psi\|} \, 
		\Bl\langle f_o(t + h_t\,\cdot) - f_o(t) + L_t h_t^k, \psi \Br\rangle 
		+ \Gamma((a+b)h_t) + \kappa_\alpha \\
	& \le & - n^{1/2} L_t h_t^{k+1/2} \|\psi\| + \Gamma((a+b)h_t) + \kappa_\alpha ,
\eea
where the latter inequality follows from Lemma~\ref{Nice Psis}~(c). Specifically let
$$
	h_t \ := \ c w_\epsilon(t)^2 \rho_n^{1/k}
$$
for some positive constant $c$ to be specified later. By continuity of $\Dk f_o$, the weight function $w_\epsilon$ is bounded away from zero and infinity. Hence $h_t\to 0$ and $L_t \max(\Dk f_o(t),\epsilon)^{-1}\to 1$, uniformly in $t\in[0,1]$. In particular,
\bea
	\Gamma((a+b)h_t) 
	& \le & (k+1/2)^{-1/2} \log(en)^{1/2}	
		\quad\mbox{for } n \ge n_o , \\
	n^{1/2} L_t h_t^{k+1/2} \|\psi\| 
	& \ge & c^{k+1/2} \|\psi\| \log(en)^{1/2} , \\
	L_t h_t^k 
	& \le & w_\epsilon(t)^{-1} c^k (1 + b_n) \rho_n ,
\eea
where $n_o$ and $b_n$ are positive numbers depending only on $f_o$, $\epsilon$ and $c$ such that $b_n\to 0$. Consequently, for $n \ge n_o$,
$$
	ah_t \leq t \leq 1 - bh_t	\quad\mbox{and}\quad	
	(f_o - \hat\ell)(t) w_\epsilon(t) \ \ge \ c^k (1 + b_n) \rho_n
$$
implies that
$$
	{\psi W(h_t,t) \over h_t^{1/2} \|\psi\|} 
	\ \le \ - \Bl( c^{k+1/2} \|\psi\| - (k+1/2)^{-1/2} \Br) \log(en)^{1/2} 
		+ \kappa_\alpha .
$$
Whenever $c > (\Delta^{(\ell)})^{1/k}$, the right-hand side of the preceding inequality tends to minus infinity, while the random variable on the left-hand side has mean zero and variance one. Since the limit of $c^k (1 + b_n)$ can be arbitrarily close to $\Delta^{(\ell)}$, these considerations show that $(f_o - \hat\ell)(t) w_\epsilon(t) \le (\Delta^{(\ell)} + o_p(1)) \rho_n$ for any fixed $t \in \left]0,1\right[$.

If $n$ is sufficiently large, then $ah_t \le t \le 1-bh_t$ and
$$
	{\psi W(h_t,t) \over h_t^{1/2} \|\psi\|} \ \ge \ - T(-\psi) - \Gamma((a+b)h_t)
$$
for all $t\in[\epsilon,1-\epsilon]$. Consequently,
$$
	\sup_{t\in[\epsilon,1-\epsilon]} (f_o - \hat\ell)(t( w_\epsilon(t) \ \geq \ c^k (1 + b_n)
$$
implies that
\bea
	T(-\psi) 
	& \ge & n^{1/2} L_t h_t^{k+1/2} \|\psi\| - 2 \Gamma((a+b)h_t) - \kappa_\alpha \\
	& \ge & \Bl( c^{k+1/2} \|\psi\| - 2 (k+1/2)^{-1/2} \Br) \log(en)^{1/2} - \kappa_\alpha .
\eea
Whenever $c > 2^{1/(k+1/2)} (\Delta^{(\ell)})^{1/k}$, the right hand side of the preceding inequality tends to infinity. Since the limit of $c^k (1 + b_n)$ can be arbitrarily close to $2^{k/(k+1/2)} \Delta^{(\ell)}$, these considerations reveal that $\|(f_o - \hat\ell)w_\epsilon\|^+_{\epsilon,1-\epsilon}$ is not greater than $\Bl( 2^{k/(k+1/2)} \Delta^{(\ell)} + o_p(1) \Br) \rho_n$.	\hfill	$\Box$

%=============================
\section{Some decision theory}
\label{Decision Theory}
%=============================

Let $X = (X_i)_{i=1}^m$ be a random vector with distribution $\NN_m(\theta, I)$. In what follows we consider tests $\phi : \R^m \to [0,1]$ and confidence sets
$$
	S \ = \ S_1 \times S_2 \times \cdots \times S_m
$$
for $\theta$ with random intervals $S_j \subset \R$. The conditional distribution of $S$, given $X$, does not depend on $\theta$. The possibility of randomized confidence sets $S$, i.e.~confidence sets not just being a function of $X$, has to be included for technical reasons. Unless specified differently, asymptotic statements in this section refer to $m \to \infty$.

\begin{Theorem}	\label{CS Sup-Diameter}
Let $c_m := (2\log m)^{1/2}$. There are universal positive numbers $b_m$ with $b_m\to 0$ such that the following two inequalities are satisfied:

\noindent
{\bf (a)} For arbitrary tests $\phi$,
$$
	\min_{j=1,\ldots,m} \Ex_{(c_m - b_m) e_j} \phi(X) - \Ex_0 \phi(X) 
	\ \le \ b_m ,
$$
where $e_1,e_2, \ldots, e_m$ denotes the standard basis of $\R^m$.

\noindent
{\bf (b)} For arbitrary confidence sets $S$ as above,
$$
	\min_{\theta \in [-c_m,c_m]^m} \, 
		\Pr_\theta \Bl\{ \theta \in S 
			\mbox{ and } \max_{j = 1,\ldots,m} {\rm length}(S_j) < 2 (c_m - b_m) \Br\} 
	\ \le \ b_m .
$$
\end{Theorem}

\noindent
{\bf Proof of Theorem~\ref{CS Sup-Diameter}.} Part~(a) is classical and can be proved by a Bayesian argument; see for instance Ingster~(1993) or D\"umbgen and Spokoiny~(2001). In order to prove part~(b) we also consider a Bayesian model: Let $\theta$ have independent components each of which is uniformly distributed on the three-point set $K_m := \{-\kappa_m,0,\kappa_m\}$, where $\kappa_m := c_m - b_m$ with constants $b_m \in [0,c_m]$ to be specified later on. Let $\LL(X \,|\, \theta) = \NN_m(\theta, I)$. Let $\Pr(\cdot), \Ex(\cdot)$ denote probabilities and expectations in this Bayesian context, whereas $\Pr_\theta(\cdot), \Ex_\theta(\cdot)$ are used in case of a fixed parameter $\theta$. For any confidence set $S$,
\bea
	\lefteqn{ \min_{\theta \in [-c_m,c_m]^m} \, 
		\Pr_\theta \Bl\{ \theta \in S 
			\mbox{ and } \max_{j = 1,\ldots,m} {\rm length}(S_j) < 2 \kappa_m \Br\} } \\
	& \le & \Pr \Bl\{ \theta \in S 
		\mbox{ and } \max_{j = 1,\ldots,m} {\rm length}(S_j) < 2 \kappa_m \Br\} 
		\ \le \ \Pr\{\theta \in \tilde{S}\} ,
\eea
where
$$
	\tilde{S} \ := \ \left\{\barr{cl} 
		S & \mbox{if } \displaystyle\max_{j = 1,\ldots,m} {\rm length}(S_j) < 2 \kappa_m , \\
		\{0\} \times \cdots \times \{0\} & \mbox{else} .
	\earr\right.
$$
The conditional distribution of $\theta$ given $(X,S)$ is also a product of $m$ probability measures: For any $\eta\in K_m^m$,
$$
	\Pr(\theta = \eta \,|\, X,S) 
	\ = \ \prod_{i=1}^m g(\eta_i \,|\, X_i)	
	\quad\mbox{with}\quad	
	g(z \,|\, x) \ := \ { \exp(- (x - z)^2/2) \over \sum_{y \in K_m} \exp(- (x - y)^2/2) } .
$$
Since each factor $\tilde{S}_j$ of $\tilde{S}$ contains at most two points from $K_m$,
\bea
	\Pr\{\theta \in \tilde{S}\} 
	& = & \Ex \Pr(\theta \in \tilde{S} \,|\, X,S) \\
	& \le & \Ex \max_{\eta \in K_m^m} \, 
		\Pr(\theta_i \neq \eta_i \mbox{ for } i=1,\ldots,m \,|\, X,S) \\
	& = & \Ex \prod_{i=1}^m \Bl( 1 - \min_{z\in K_m} g(z \,|\, X_i) \Br) \\
	& = & \Bl( 1 - \Ex \min_{z\in K_m} g(z \,|\, X_1) \Br)^m \\
	& \le & \Bl( 1 - 3^{-1} \Ex \min_{z\in K_m} \, \exp(- (X_1 - z)^2/2) \Br)^m .
\eea
The latter expectation can be bounded from below as follows:
\bea
	\lefteqn{ 3^{-1} \Ex \min_{z\in K_m} \, \exp(- (X_1 - z)^2/2) } \\
	& \ge & 3^{-1} \Pr\{|X_1| \le b_m/2\} \exp(- (\kappa_m + b_m/2)^2/2) \\
	& \ge & 3^{-1} \Pr\{|\theta_1| = 0, |X_1| \le b_m/2\} \exp(- (c_m - b_m/2)^2/2) \\
	& = & 9^{-1} (2\pi)^{-1/2} (b_m + O(b_m^2)) \exp(c_m b_m/2 - b_m^2/8) m^{-1} .
\eea
In case of $b_m := 1\{m > 1\} c_m^{-1/2} = o(1)$ the latter bound is easily seen to be $a_m m^{-1}$ with $a_m = a_m(b_m)\to\infty$. Thus
$$
	\Pr\{\theta \in \tilde{S}\}	
	\ \le \ (1 - a_m m^{-1})^m 
	\ \to \ 0 .
$$
Replacing $b_m$ with $\max\{b_m, (1 - a_m m^{-1})^m\}$ yields the assertion of part~(b).	\hfill	$\Box$

%======================================
\section{Related optimization problems}
\label{Optimization}
%======================================

As in Section~\ref{GG smooth} let $(\GG,k)$ be either $(\GGiso,1)$ or $(\GGconv,2)$. In view of future applications to other regression models we extend our framework slightly and consider $\langle g,h\rangle := \int gh \, d\mu$, $\|g\| := \langle g,g\rangle^{1/2}$ for some measure $\mu$ on the real line such that $\mu(C)<\infty$ for bounded intervals $C\subset\R$.

Let $\psi$ be some bounded function on the real line with $\psi(x) = 0$ for $x \not\in [-a,b]$ and $\langle 1,\psi\rangle \ge 0$, where $a,b \ge 0$. The next lemma provides sufficient conditions for one of the following two requirements:
\bean
	\langle g, \psi\rangle & \le & g(0) \langle 1,\psi\rangle	
		\quad\mbox{whenever } g \in \GG , 1_{[-a,b]} g \in L^1(\mu) , 
	\label{OrthoL} \\
	\langle g, \psi\rangle & \ge & g(0) \langle 1,\psi\rangle	
		\quad\mbox{whenever } g \in \GG , 1_{[-a,b]} g \in L^1(\mu) .
	\label{OrthoU}
\eean

\begin{Lemma}	\label{Ortho}
Let $\GG = \GGiso$ and $\psi \ge 0$. Then $b = 0$ entails condition~(\ref{OrthoL}), while $a = 0$ implies condition~(\ref{OrthoU}).

Let $\GG = \GGconv$ and $\int_{-\infty}^\infty x \psi(x) \, \mu(dx) = 0$. Condition~(\ref{OrthoU}) is satisfied if $\psi\ge 0$. On the other hand, condition~(\ref{OrthoL}) is a consequence of the following two requirements: 
$\int x^\pm\psi(x) \, \mu(dx) = 0$ and
$$
	\psi \ \left\{\barr{cl} 
		\ge 0 & \mbox{on } [c,d] \\
		\le 0 & \mbox{on } \R\setminus[c,d]
	\earr\right.
$$
for some numbers $c < 0 < d$, where $\mu([-a,c]), \mu([d,b]) > 0$. (Here $y^+ := \max(y,0)$ and $y^- := \max(-y,0)$.)
\end{Lemma}

With Lemma~\ref{Ortho} at hand one can solve two mimimization problems leading to the special kernels in (\ref{PsiL}) and (\ref{PsiU}). In both cases we consider two disjoint convex sets $\GG_o,\GG_A \subset \GG$ and construct functions $G_o \in \GG_o$, $G_A \in \GG_A$ such that
\be
	\|G_o - G_A\| \ = \ \min_{g_o\in\GG_o, \, g_A\in\GG_A} \, \|g_o-g_A\| .
	\label{OptimizationG}
\ee

\begin{Theorem}	\label{OptimizationL}
Let $\GG_o := \Bl\{ g \in \GG : g(0) \le -1 \Br\}$ and $\GG_A := \Bl\{ g \in \GG \cap \HH_{k,1} : g(0) \ge 0 \Br\}$. In case of $\GG=\GGiso$ let $G_A(x) := x$ and
$$
	G_o(x) \ := \ \left\{\barr{cl} 
		-1     & \mbox{if } x \in [-1,0] , \\
		G_A(x) & \mbox{else} .
	\earr\right.
$$
In case of $\GG=\GGconv$ let $G_A(x) := x^2/2$ and
$$
	G_o(x) \ := \ \left\{\barr{cl} 
		-1 + (a/2+1/a) x^- + (b/2+1/b) x^+ & \mbox{if } x \in [-a,b] , \\
		G_A(x) & \mbox{else} ,
	\earr\right.
$$
where $a,b \ge 2^{1/2}$ are chosen such that $\int x^\pm(G_A-G_o)(x)\,\mu(dx) = 0$.

Then equation~(\ref{OptimizationG}) holds in both cases. More precisely, the function $\psi := G_A - G_o$ satisfies the inequalities $\langle 1,\psi\rangle \ge \|\psi\|^2$, (\ref{OrthoL}) and
\be
	\langle g,\psi\rangle \ \ge \ \|\psi\|^2 - \langle 1, \psi\rangle
	\quad\mbox{whenever } g \in \HH_{k,1}, g(0) \ge 0 .
	\label{OrthoL2}
\ee
\end{Theorem}

In case of $\mu$ being Lebesgue measure, $\psi = G_A-G_o$ coincides with the function $\psi^{(\ell)}$ in (\ref{PsiL}), where $a = b = 2$.

\begin{Theorem}	\label{OptimizationU}
Let $\GG_o := \Bl\{ g \in \GG : g(0) \ge 1 \Br\}$, $\GG_A := \Bl\{ g \in \GG \cap \HH_{k,1} : g(0) \le 0 \Br\}$, and define $G_A$ as in Theorem~\ref{OptimizationL}. In case of $\GG=\GGiso$ let
$$
	G_o(x) \ := \ \left\{\barr{cl} 
		0 & \mbox{if } x \in [0,1] , \\
		G_A(x) & \mbox{else} .
	\earr\right.
$$
In case of $\GG=\GGconv$ suppose that $\mu(\left]-\infty,0\right[),\mu(\left]0,\infty\right[)>0$ and let
$$
	G_o(x) \ := \ \left\{\barr{cl} 
		1 + cx & \mbox{if } x \in [-a,b] , \\
		G_A(x) & \mbox{else} ,
	\earr\right.
$$
where $a := -c + (c^2 + 2)^{1/2}$, $ b := c + (c^2 + 2)^{1/2}$, and $c$ is chosen such that $\int x (G_o - G_A)(x) \, \mu(dx) = 0$.

Then equation~(\ref{OptimizationG}) is satisfied in both cases. More precisely, the function $\psi := G_o-G_A$ satisfies the inequalities $\langle 1,\psi\rangle \ge \|\psi\|^2$, (\ref{OrthoU}) and
\be
	\langle g,\psi\rangle \ \le \ \langle 1,\psi\rangle - \|\psi\|^2	
	\quad\mbox{whenever } g \in \HH_{k,1}, g(0) \ge 0 .
	\label{OrthoU2}
\ee
\end{Theorem}

In case of $\mu$ being Lebesgue measure, $\psi = G_o-G_A$ coincides with the function $\psi^{(u)}$ in (\ref{PsiU}), where $c = 0$ and $a = b = 2^{1/2}$.

The following lemma summarizes essential properties of the optimal kernels $\psi^{(\ell)}$ and $\psi^{(u)}$.

\begin{Lemma}	\label{Nice Psis}
Let $\psi^{(\ell)}$ and $\psi^{(u)}$ be the kernel functions in (\ref{PsiL}) and (\ref{PsiU}), and let $h,L > 0$ and $t \in \R$.

{\bf(a)} If $\GG = \GGiso$, then $\langle 1,\psi^{(\ell)}\rangle = \langle 1,\psi^{(u)}\rangle = 1/2$ and $\|\psi^{(\ell)}\|^2 = \|\psi^{(u)}\|^2 = 1/3$. If $f:\R\to\R$ satisfies $f(y) - f(x) \ge L(y-x)$ for all $x < y$, then
$$
	f - Lh^{-1} \psi^{(\ell)}_{h,t}, f + Lh^{-1} \psi^{(u)}_{h,t} \ \in \ \GGiso .
$$

{\bf(b)} If $\GG = \GGconv$, then $\langle 1,\psi^{(\ell)}\rangle = 2/3$, $\|\psi^{(\ell)}\|^2 = 8/15$, $\langle 1,\psi^{(u)}\rangle = 2^{2.5}/3$ and $\|\psi^{(u)}\|^2 = 2^{4.5}/15$. Let $f:\R\to\R$ be absolutely continuous with derivative $f'$ such that $f'(y)-f'(x)\ge L(y-x)$ for all $x<y$. Then
$$
	f - Lh^{-2} \psi^{(\ell)}_{h,t}, f + Lh^{-2} \psi^{(u)}_{h,t} 
	\ \in \ \GGconv .
$$

{\bf(c)} In general, for any function $f \in \HH_{k,L}$,
\bea
	\left\langle f(t + h\,\cdot) - r + Lh^k, \psi^{(\ell)} \right\rangle 
	& \ge & Lh^k \|\psi^{(\ell)}\|^2	\quad\mbox{if } f(t) \ge r , \\	
	\left\langle f(t + h\,\cdot) - r - Lh^k, \psi^{(u)} \right\rangle 
	& \le & -Lh^k \|\psi^{(u)}\|^2	\quad\mbox{if } f(t) \le r .
\eea
\end{Lemma}

\noindent
{\bf Proof of Lemma~\ref{Ortho}.} The assertions for $\GG = \GGiso$ are a simple consequence of $g \le g(0)$ on $\left]-\infty,0\right]$ and $g \ge g(0)$ on $\left[0,\infty\right[$.

Now let $\GG = \GGconv$. If $\psi \ge 0$ and $\int x \psi(x) \, \mu(dx) = 0$, then Condition~(\ref{OrthoU}) follows from Jensen's inequality applied to the probability measure $P(dx) = \langle 1,\psi\rangle^{-1} \psi(x) \, \mu(dx)$.

On the other hand, suppose that $\psi \ge 0$ on $[c,d]$ and $\psi \le 0$ on $\R\setminus [c,d]$, where $c < 0 < d$ and $\mu([-a,c]), \mu([d,b]) > 0$. For $g \in \GGconv$ with $1_[-a,b] g \in L^1(\mu)$, both $g(c)$ and $g(d)$ have to be finite, and we define
$$
	\tilde{g}(x) \ := \ 
		g(x) - \left\{\barr{cl} 
			d^{-1} (g(d)-g(0)) x & \mbox{if } x \ge 0 , \\
			c^{-1} (g(c)-g(0)) x & \mbox{if } x \le 0 .
		\earr\right.
$$
By convexity of $g$, this auxiliary function $\tilde{g}$ satisfies $\tilde{g} \le g(0)$ on $[c,d]$ and $\tilde{g} \ge g(0)$ on $\R \setminus [c,d]$. Thus $\langle \tilde{g}, \psi\rangle \le g(0) \langle 1,\psi\rangle$. If in addition $\int x^\pm \psi(x) \, \mu(dx) = 0$, then $\langle g,\psi\rangle = \langle\tilde{g},\psi\rangle$.	\hfill	$\Box$

\noindent
{\bf{Proof of Theorem~\ref{OptimizationL}.}} One can easily deduce from Lemma~\ref{Ortho} that the function $\psi = G_A - G_o$ satisfies inequality~(\ref{OrthoL}). But $G_A$ is an extremal point of $\GG_A$ in the sense that
$$
	G_A - g \ \in \ \GG	\quad\mbox{for any } g \in \HH_{k,1} .
$$
For let $x < y$. If $\GG = \GGiso$, then
$$
	(G_A-g)(y) - (G_A-g)(x) 
	\ = \ y - x - (g(y) - g(x)) \ \ge \ y - x - |y-x| \ = \ 0 ,
$$
whence $G_A - g$ is non-decreasing. In case of $\GG = \GGconv$ the same argument applies to the first derivative of $G_A - g$. Together with (\ref{OrthoL}) this implies that
\bea
	\langle g,\psi\rangle 
	& = & \langle G_A, \psi\rangle - \langle G_A-g,\psi\rangle \\
	& \ge & \langle G_A, \psi\rangle - (G_A - g)(0) \langle 1,\psi\rangle \\
	& = & \langle G_A, \psi\rangle + g(0) \langle 1,\psi\rangle \\
	& = & \|\psi\|^2 + \langle G_o,\psi\rangle + g(0) \langle 1,\psi\rangle \\
	& = & \|\psi\|^2 + (g(0)-1) \langle 1,\psi\rangle .
\eea
The latter equation follows from $\langle G_o,\psi\rangle = \langle -1,\psi\rangle$, which is easily verified. The special case $g = 0$ yields the inequality $\langle 1,\psi\rangle \geq \|\psi\|^2$. Then inequality (\ref{OrthoL2}) becomes obvious.

It remains to be shown that in case of $\GG = \GGconv$ there exist numbers $a,b \ge 2^{1/2}$ such that $\psi = \psi(\cdot,a,b)$ satisfies $\int x^\pm \psi(x) \, \mu(dx) = 0$. In fact, for any fixed $x$ the number $\psi(x,a,b) \le 1$ can be shown to be continuous and decreasing in $a$ and $b$. Precisely, $\psi(0,a,b) = 1$ and $\lim_{a\to\infty} \psi(x,a,\cdot) = \lim_{b\to\infty} \psi(y,\cdot,b) = -\infty$ for $x < 0 < y$. Hence the assertion is a consequence of monotone convergence.	\hfill	$\Box$

\noindent
{\bf Proof of Theorem~\ref{OptimizationU}.} This proof is analogous to the proof of Theorem~\ref{OptimizationL} and thus omitted.	\hfill	$\Box$

\noindent
{\bf Proof of Lemma~\ref{Nice Psis}.} The calculations of $\langle 1,\psi\rangle$ and $\|\psi\|^2$ are elementary and thus omitted. Elementary calculations show that $g := - Lh^{-k}\psi^{(\ell)}_{t,h}$ as well as $g := Lh^{-k}\psi^{(u)}$ satisfies
$$
	\left.\barr{c} 
		g(y) - g(x) \\
		g'(y)-g'(x) 
	\earr\right\} \ \ge - L(y-x)	
	\quad\mbox{if } \GG = \left\{\barr{c} 
		\GGiso , \\
		\GGconv , 
	\earr\right.
$$
where $g'(x)$ denotes any number between the right- and left-sided derivative of $g$ at $x$. Thus $f+g$ belongs to $\GG$, whenever $f$ satisfies the inequalities stated in parts~(a) and (b).

As for part~(c), for $f \in \HH_{k,L}$ and $t\in \R$, $h,c > 0$ the function $c f(t + h\,\cdot)$ belongs to $\HH_{k,cLh^k}$. If we take $c := (Lh^k)^{-1}$, the inequality (\ref{OrthoL2}) implies that
\bea
	\left\langle f(t + h\,\cdot) - r + Lh^k,\psi^{(\ell)} \right\rangle 
	& = & Lh^k \left\langle c(f(t + h\,\cdot)-f(t)) + 1, \psi^{(\ell)} \right\rangle \\
	& \ge & Lh^k \|\psi^{(\ell)}\|^2 .
\eea
Analogously one can deduce the lower bound for $\left\langle f(t+h\,\cdot)-r-Lh^k,\psi^{(u)}\right\rangle$.	\hfill	$\Box$

\bigskip

{\bf Acknowledgements.} The author is grateful to Lars H\"omke for his assistance in Section~\ref{Examples}. Constructive comments of a referee and an associate editor helped to improve the presentation. This work has been supported by Deutsche Forschungsgemeinschaft, grant Du\,238/5-1.

%========================
\subsection*{References}
%========================

\begin{description}
\item[]{\sc Bickel, P.J. and M. Rosenblatt} (1973).
		On some global measures of the deviations of density function estimates. \
		{\sl Ann. Statist.}~{\bf 1}, 1071--1095
\item[]{\sc Brown, L.D. and M.G. Low} (1996). 
		Asymptotic equivalence of nonparametric regression and white noise. \ 
		{\sl Ann. Statist.}~{\bf 24}, 2384--2398.
\item[]{\sc Davies, P.L. (1995).} 
		Data features. \ 
		{\sl Statistica Neerlandica}~{\bf 49}, 185--245
\item[]{\sc Donoho, D.L. (1988).} 
		One-sided inference about functionals of a density. \ 
		{\sl Ann. Statist.}~{\bf 16}, 1390--1420
\item[]{\sc D\"umbgen, L.} (1998). 
		New goodness-of-fit tests and their application to nonparametric confidence sets. \ 
		{\sl Ann. Statist.}~{\bf 26}, 288--314
\item[]{\sc D\"umbgen, L.} (2007).
		Confidence bands for convex median functions using sign tests.\\
		In: \textsl{Asymptotics: Particles, Processes and Inverse Problems
		(E. Cator, G. Jongbloed, C. Kraai\-kamp, R. Lopuha{\"a}, J.A. Wellner, eds.)},
		pp. 85-100. \
		Lecture Notes - Monograph Series \textbf{55}, IMS, Hayward, USA.
\item[]{\sc D\"umbgen, L. and R.B. Johns} (2004).
		Confidence bands for isotonic median functions using sign tests. \
		{\sl J. Comp. Graph. Statist.}~{\bf 13}, 519--533
\item[]{\sc D\"umbgen, L. and V.G. Spokoiny} (2001). 
		Multiscale testing of qualitative hypotheses. \ 
		{\sl Ann. Statist.}~{\bf 29}, 124--152
\item[]{\sc Eubank, R.L. and P.L. Speckman} (1993).
		Confidence bands in nonparametric regression. \
		{\sl J. Amer. Statist. Assoc.}~{\bf 88}, 1287--1301
\item[]{\sc Fan, J. and W. Zhang} (2000).
		Simultaneous confidence bands and hypothesis testing in varying-coefficient models. \
		{\sl Scand. J. Statist.}~{\bf 27}, 715--731
\item[]{\sc Fan, J., C. Zhang and J. Zhang} (2001).
		Generalized likelihood ratio statistics and Wilks phenomenon. \
		{\sl Ann. Statist.}~{\bf 29}, 153--193
\item[]{\sc Grama, I. and M. Nussbaum} (1998). 
		Asymptotic equivalence for nonparametric generalized linear models. \ 
		{\sl Prob. Theory and Related Fields}~{\bf 111}, 167--214.
\item[]{\sc H\"ardle, W. and J.S. Marron} (1991).
		Bootstrap simultaneous error bars for nonparametric regression. \
		{\sl Ann. Statist.}~{\bf 19}, 778--796
\item[]{\sc Hall, P. and D.M. Titterington} (1988).
		On confidence bands in nonparametric density estimation. \
		{\sl J. Multivar. Anal.}~{\bf 27}, 228--254
\item[]{\sc Hart, J.D.} (1997).
		{\sl Nonparametric Smoothing and Lack-of-Fit Tests.} \
		Springer, New York
\item[]{\sc Hengartner, N.W. and P.B. Stark (1995).} 
		Finite-sample confidence envelopes for shape-restricted densities. \ 
		{\sl Ann. Statist.}~{\bf 23}, 525--550
\item[]{\sc Khas'minskii, R.Z.} (1978). 
		A lower bound on the risks of nonparametric estimates of densities 
		in the uniform metric. \ 
		{\sl Theory Prob. Appl.}~{\bf 23}, 794--798
\item[]{\sc Knafl, G., J. Sachs and D. Ylvisaker} (1985).
		Confidence bands for regression functions. \
		{\sl J. Amer. Statist. Assoc.}~{\bf 80}, 683--691
\item[]{\sc Nussbaum, M.} (1996). 
		Asymptotic equivalence of density estimation and white noise. \ 
		{\sl Ann. Statist.}~{\bf 24}, 2399--2430.
\item[]{\sc Robertson, T., F.T. Wright and R.L. Dykstra} (1988). 
		{\sl Order Restricted Statistical Inference.} \ 
		Wiley, New York
\item[]{\sc Ingster, Y.I. (1993).} 
		Asymptotically minimax hypothesis testing for nonparametric alternatives, I-III. \ 
		{\sl Math. Methods Statist.}~{\bf 2}; 85-114, 171--189, 249-268.
\end{description}

%=============
\end{document}